\documentclass[preprint]{imsart}  

\usepackage{enumerate}
\usepackage{amsmath,latexsym,amssymb,amsthm}
\usepackage{graphicx}	
\usepackage{graphics}											
\usepackage{import}												
\usepackage{graphicx}												
\usepackage{color}													
\usepackage{multirow}
\usepackage{chngpage}
\RequirePackage[colorlinks,citecolor=blue,urlcolor=blue]{hyperref}

\arxiv{arXiv:1609.04481}

\newcommand{\rmd}{{\rm d}}
\newcommand{\rmi}{{\rm i}}

\newcommand{\eqd}{\stackrel{D}{=}}

\newcommand{\halmos}{\quad\hfill\mbox{$\Box$}}

\newcommand{\EE}{\mathbb{E}}
\newcommand{\DD}{\mathbb{D}}
\newcommand{\RR}{\mathbb{R}}

\newcommand{\ZZ}{\mathbb{Z}}
\newcommand{\CC}{\mathbb{C}}

\newcommand{\PP}{\mathbb{P}}

\newcommand{\mySS}{\mathbb{S}}

\newcommand{\DDD}{{\cal D}}

\newcommand{\FFF}{{\cal F}}
\newcommand{\GGG}{{\cal G}}

\newcommand{\KKK}{{\cal K}}

\newcommand{\NNN}{{\cal N}}

\newcommand{\RRR}{{\cal R}}
\newcommand{\SSS}{{\cal S}}
\newcommand{\TTT}{{\cal T}}
\newcommand{\UUU}{{\cal U}}
\newcommand{\VVV}{{\cal V}}

\newcommand{\XXX}{{\cal X}}
\newcommand{\YYY}{{\cal Y}}
\newcommand{\ZZZ}{{\cal Z}}
\newcommand{\KKKK}{\mathfrak{K}}
\newcommand{\GGC}{GGC}
\newcommand{\VGGC}{VGG}
\newcommand{\Levy}{L\'evy }
\newcommand{\skal}[2]{\left\langle #1,#2\right\rangle}

\newcommand{\eins}{{\bf 1}}

\newcommand{\bfnull}{{\bf 0}}

\newcommand{\bfc}{{\bf c}}
\newcommand{\bfe}{{\bf e}}
\newcommand{\bfm}{{\bf m}}
\newcommand{\bfs}{{\bf s}}
\newcommand{\bft}{{\bf t}}
\newcommand{\bfx}{{\bf x}}
\newcommand{\bfy}{{\bf y}}

\newcommand{\bfz}{{\bf z}}
\newcommand{\bfd}{{\bf d}}
\newcommand{\bfu}{{\bf u}}
\newcommand{\bfv}{{\bf v}}
\newcommand{\bfalpha}{\boldsymbol{\alpha}}
\newcommand{\bfdelta}{\boldsymbol{\delta}}

\newcommand{\bflambda}{\boldsymbol{\lambda}}
\newcommand{\bfmu}{\boldsymbol{\mu}}
\newcommand{\bfpi}{\boldsymbol{\pi}}

\newcommand{\bftheta}{\boldsymbol{\theta}}

\newcommand{\myCov}{{\rm Cov}}
\newcommand{\myVar}{{\rm Var}}

\newcommand{\diag}{\mbox{diag}}

\newcommand{\spur}{\mbox{trace}}

\newcommand{\tr}{\diamond}
\newcommand{\su}{\odot}

\numberwithin{equation}{section}

\newtheorem{theorem}{Theorem}
\newtheorem{lemma}{Lemma}
\newtheorem{remark0}{\bf Remark}
\newenvironment{remark}{\begin{remark0}\em}{\end{remark0}\par}

\newtheorem{definition}{Definition}
\newtheorem{proposition}{Proposition}
\numberwithin{theorem}{section}
\numberwithin{proposition}{section}
\numberwithin{lemma}{section}
\numberwithin{corollary}{section}
\numberwithin{remark0}{section}
\numberwithin{definition}{section}

\begin{document}

\begin{frontmatter}

\title{Weak Subordination of Multivariate L\'evy Processes and Variance Generalised Gamma Convolutions}

\runtitle{Weak Subordination of Multivariate L\'evy Processes}
\begin{aug}
\author{\fnms{Boris} \snm{Buchmann}\thanksref{a}\ead[label=e1]{Boris.Buchmann@anu.edu.au}}
\and \author{\fnms{Kevin~W.} \snm{Lu}\thanksref{b}\ead[label=e2]{kevin.lu@anu.edu.au}}
\and
\author{\fnms{Dilip~B.} \snm{Madan}\thanksref{cb}\ead[label=e3]{dbm@rhsmith.umd.edu}}

\address[a]{Research School of Finance, Actuarial Studies \& Statistics,
Australian National University, ACT 0200, Australia.
\printead{e1}}
\address[b]{Mathematical Sciences Institute, Australian National University, ACT 0200, Australia.
\printead{e2}}
\address[cb]{Robert H. Smith School of Business, University of Maryland, College Park, MD. 20742, USA.
\printead{e3}}

\runauthor{B.Buchmann et al.}
\end{aug}

\begin{abstract}
	Subordinating a multivariate L\'evy process, the subordinate, with a univariate subordinator gives rise to a pathwise construction of a new L\'evy process, provided
	the subordinator and the subordinate are independent processes.
	The variance-gamma model in finance was generated accordingly from a Brownian motion and a gamma process.
	Alternatively, multivariate subordination can be used to create L\'evy processes, but this requires the subordinate to have independent components.
	In this paper, we show that there exists another operation acting on pairs $(T,X)$ of L\'evy processes which creates a L\'evy process $X\odot T$.
	Here, $T$ is a subordinator, but $X$ is an arbitrary L\'evy process with possibly dependent components.
	We show that this method is an extension of both univariate and multivariate subordination and provide two applications. We illustrate our methods giving a weak formulation of the variance-$\bfalpha$-gamma process that exhibits a wider range of dependence than using traditional subordination.
	Also, the variance generalised gamma convolution class of L\'evy processes formed by subordinating Brownian motion with Thorin subordinators is further extended using weak subordination.
\end{abstract}



\end{frontmatter}


\section{Introduction}\label{secintro}
The subordination of L\'evy processes has many important applications.
In mathematical finance, for instance, it acts as a time change that models the flow of information, measuring time in volume of trade as opposed to real time. This idea was initiated by~\cite{MaSe90} who introduced the variance-gamma process for modelling stock prices, where the subordinate is Brownian motion and the subordinator is a gamma process. Multivariate subordination can be applied to model dependence across multivariate L\'evy processes, where the components may have common and/or idiosyncratic time changes. We refer the reader to~\cite{BKMS16} for a thorough discussion of traditional subordination and its applications.

Let $T=(T_1,\dots,T_n)$ be an $n$-dimensional subordinator, and $X=(X_1,\dots,X_n)$ be another $n$-dimensional L\'evy process called the subordinate. Subordination is the operation that produces the $n$-dimensional process $X\circ T$ defined by
\[(X\circ T)(t):=(X_1(T_1(t)),\dots,X_n(T_n(t)))\,,\quad t\ge 0\,.\]
If $T$ and $X$ are {\em independent}, then there are two important special cases where $X\circ T$ is again a L\'evy process:\\[1mm]
\noindent$\bullet$~Traditional Bochner subordination, where $T_1=T_2=\dots=T_n$ are {\em indistinguishable}~\cite{BS10,R065,s,Zo58}.\\
\noindent$\bullet$~$T$ is multivariate but $X$ has {\em independent components}~$X_1,\dots,X_n$~\cite{BPS01}.

Thus, for strictly multivariate subordination, that is $n\ge 2$, while $T$ does not have indistinguishable components, we have to restrict the class of admissible subordinates $X$ to L\'evy processes with {\em independent components}, which is,
as we show in Proposition~\ref{propindepindist} below, in some cases necessary if we are to stay in the class of L\'evy processes.

In the present paper, we show that there exists an operation that extends the traditional notion of subordination by assigning the distribution of a L\'evy process $X\odot T$ to the pair $(T,X)$ of L\'evy processes. The weakly subordinated process $X\odot T$ is a general L\'evy process, it inherits jumps from the multivariate subordinator $T$, which resembles subordination, and our new operation reduces to subordination when the components of $X$ are independent or the components of $T$ are indistinguishable.

The remaining parts of the paper are organised as follows. In Section~\ref{secmain}, we define weak subordination, show its existence and that there is a stronger pathwise interpretation, based on marked point processes of jumps, which we call semi-strong subordination. We review some properties of gamma and variance-gamma processes. Further, we introduce a weakly subordinated version of the variance-$\bfalpha$-gamma process as an extension of the
strongly subordinated version in \cite{Se08}. We develop this new class in a number of remarks throughout the paper to illustrate
our machinery.

Section~\ref{secprop} contains a number of results concerning the relation between traditional subordination and weak subordination. In particular, we show that weak subordination extends traditional subordination and is consistent with projecting to marginal distributions, like traditional subordination. However, there are also differences between both notions. To highlight these, we provide formulae for the first and second moments and covariances for weakly subordinated processes.

In Section~\ref{secVGGC}, we exemplify the unifying nature of weak subordination as illustrated using variance  generalised gamma convolutions. In~\cite{Gr07}, a class of processes was introduced by subordinating $n$-dimensional Brownian motion with univariate subordinators taken from Thorin's~\cite{Th77a,Th77b} class of {\em generalised gamma convolutions}~($GGC$). This class has been coined $VGG^{n,1}$ in~\cite{BKMS16} who complemented it with their $VGG^{n,n}$-class, obtained by subordinating $n$-dimensional Brownian motion with independent components with $n$-dimen\-sionional $GGC$-subordinators. Using weak subordination, we introduce a {\em weak $VGG^n$}-class of \Levy processes as a natural superclass of the $VGG^{n,1}$ and $VGG^{n,n}$-classes. Unifying the results in~\cite{BKMS16}, we provide formulae for the associated characteristic function and L\'evy measure.

Section~\ref{secproofs} contains technical proofs.
\section{Main Results}\label{secmain}
Let $\RR^n$ be $n$-dimensional Euclidean space whose elements are row vectors $\bfx=(x_1,\dots,x_n)$, with canonical basis $\{\bfe_k:1\le k\le n\}$, and let
$\bfe:=(1,\dots,1)\in\RR^n$.
Let $\bfx',\Sigma'$ denote the transpose of a vector $\bfx$ and a matrix $\Sigma$, respectively.
Let $\skal\bfx\bfy\!=\!\bfx\bfy'$ denote the Euclidean product with Euclidean norm $\|\bfx\|^2\!=\!\skal \bfx\bfx=\bfx\bfx'$, and set
$\skal \bfx \bfy_\Sigma\!:=\!\bfx \Sigma\bfy'$ and
$\|\bfx\|^2_\Sigma:=\skal{\bfx}{\bfx}_\Sigma$ for $\bfx,\bfy\in\RR^n$ and $\Sigma\in\RR^{n\times n}$. Let $\DD\!:=\!\{\bfx\!\in\!\RR^n:\!\|\bfx\|\!\le\!1\}$ be the Euclidean unit ball centred at the origin.
If $A\!\subseteq\!\RR^n$, set $A_*\!:=\!A\backslash\{{\bf 0}\}$ and let $\eins_A(\omega)=\bfdelta_\omega(A)$ denote the indicator function and the Dirac measure, respectively. Further, $I:[0,\infty)^n\to [0,\infty)^n$ and $\ln\!:\!\CC\backslash(-\infty,0]\!\to\!\CC$ denote the identity function and the principal branch of the logarithm, respectively.
The decomposition of an extended real number $x\in[-\infty,\infty]$ into its positive and negative parts is denoted by $x=x^+-x^-$, where $x^+=x\vee0$ and $x^-=(-x)^+=-(x\wedge0)$.

If $\emptyset\!\neq\!J\!\subseteq\!\{1,\dots,n\}$, introduce the associated projection $\bfpi_J\!:\!\RR^n\!\to\!\RR^n$ by
$\bfpi_J(\bfx)\!:=\!\bfx\bfpi_J\!:=\!\sum_{j\in J}x_j\bfe_j$. If $\VVV$ is a Borel measure on $\RR^n$, then so is the image (push forward) measure $\VVV_J:=\VVV\circ \bfpi_J^{-1}$. If $\XXX$ is a Borel measure on $\RR^n_*$, then so is $\XXX_J$, which is constructed in the usual way:
first extend $\XXX$ to a Borel measure $\VVV$ on $\RR^n$ by setting $\VVV(\{{\bf 0}\}):=0$, second let $\XXX_J$ be the restriction of $\VVV_J$ to $\RR^n_*$.
If $J\!=\!\emptyset$, we employ the conventions $\bfpi_\emptyset\!\equiv\!{\bf 0}$, $\VVV_\emptyset\equiv0$ and $\XXX_\emptyset\equiv 0$.

The reader is referred to the monographs~\cite{Ap09,b,s} for necessary material on L\'evy pro\-cesses.
Particularly, the law of a L\'evy process $X=(X_1,\dots,X_n)=(X(t))_{t\ge 0}$  is determined by its characteristic function \[\Phi_{X(t)}(\bftheta)\,:=\,\EE\exp\{\rmi\skal\bftheta{X(t)}\}\,=\,\exp\{t\allowbreak\Psi_X(\bftheta)\}\,,\quad t\ge 0\,,\]
with L\'evy exponent $\Psi_X=\Psi$ where, $\bftheta\in\RR^n$,
\begin{equation}\label{0.1}
\Psi(\bftheta)\,=\,
\rmi \skal {\bfmu}\bftheta\!-\!\frac 12\;\|\bftheta\|^2_{\Sigma}
+\int_{\RR_*^n}\left(e^{\rmi\skal\bftheta \bfx}\!-\!1\!-\!\rmi\skal\bftheta \bfx \eins_\DD(\bfx)\right)\,\XXX(\rmd \bfx)\,.
\end{equation}
Here, $\bfmu\!=\!(\mu,\dots,\mu_n)\!\in\!\RR^n$ is a row vector, $\Sigma\!=\!(\Sigma_{kl})\!\in\!\RR^{n\times n}$ is a covariance matrix, and $\XXX$ is a nonnegative Borel measure on $\RR^n_*$ satisfying
\begin{equation}\label{Piintegrab}
\int_{\RR^n_*}\,\|\bfx\|^2\wedge 1\;\XXX(\rmd \bfx)\;<\;\infty\,.
\end{equation}
We write $X\!\sim\!L^n(\bfmu,\Sigma,\XXX)$, provided $X$ is an $n$-dimensional L\'evy process with canonical triplet
$(\bfmu,\Sigma,\XXX)$. Throughout, $B=(B_1,\dots,B_n)\sim BM^n(\bfmu,\Sigma)\allowbreak:=L^n(\bfmu,\Sigma,0)$ refers to an $n$-dimensional Brownian motion $B$ with linear drift $\EE[B(t)]=\bfmu t$ and covariance matrix Cov$(B(t))=t\Sigma$, $t\geq0$.

We write $X\sim FV^n(\bfd,\XXX)$ with $\bfd:=\bfmu\!-\!\int_{\DD_*} \bfx\,\XXX(\rmd \bfx)\;\in\;\RR^n$ denoting the drift of $X$, provided the paths of $X$ are of (locally) {\em finite variation}, equivalently, $\Sigma=0$ and
\begin{equation}\label{PiFVintegrab}
\int_{\DD_*}\|\bfx\|\,\XXX(\rmd \bfx)\quad<\quad\infty\,.
\end{equation}
Particularly, $T=(T_1,\dots,T_n)\sim S^n(\bfd,\TTT)$ refers to an $n$-dimen\-sional subordinator, that is a L\'evy process with nondecreasing components with drift $\bfd\!\in\![0,\infty)^n$.

Next, we revise some properties of gamma and variance-gamma processes. Notation is borrowed from~\cite{BKMS16}.\\[1mm]
\noindent{\bf Gamma subordinator.}~If $a,b>0$, then a subordinator $G$ is a {\em gamma subordinator} if and only if its marginal $G(t)\!\sim\!\Gamma(at,b)$, $t\!\geq\!0$
is gamma distributed with shape parameter $at$ and rate parameter $b$. A
drift-less subordinator $G$ with L\'evy measure $\GGG_{a,b}$ is a {\em gamma} subordinator with parameters $a,b$, provided
its L\'evy measure satisfies  $\GGG_{a,b}(\rmd g)\!=\!\eins_{(0,\infty)}(g)a e^{-b g}\rmd g/g$, in short,
$G\!\sim\!\Gamma_S(a,b)\!=\!S^1(0,\GGG_{a,b})$. If $G\sim\Gamma_S(a,b)$ and $\lambda\!>\!-b$, its Laplace exponent is $-\ln\EE[\exp\{-\lambda G(t)\}]=at\ln\{(b\!+\!\lambda)/b\}$.

If $a\!=\!b$, we refer to $G$ as a {\em standard} gamma subordinator, in short, $G\!\sim\!\Gamma_S(b):=\Gamma_S(b,b)$
and its L\'evy measure is denoted by $\GGG_b$. A gamma subordinator $G$ is a standard gamma subordinator if and only if $\EE[G(1)]=1$.\\[1mm]
\noindent{\bf Variance-gamma process.}~Let $b\!>\!0$, $\bfmu\!\in\!\RR^n$ and $\Sigma\!\in\!\RR^{n\times n}$ be a covariance matrix. For a Brownian motion $B\sim BM^n(\bfmu,\Sigma)$ independent of a gamma subordinator $G\sim \Gamma_S(b)$, we call $V$ a {\em variance-gamma $(VG^n)$ process}~\cite{MaSe90} with parameters $b,\bfmu,\Sigma$, if
\begin{equation*}\label{defVG} V\eqd B\circ (G\bfe)\sim VG^n(b,\bfmu,\Sigma)=BM^n(\bfmu,\Sigma)\circ(\Gamma_S(b)\bfe)\,.
\end{equation*}
An $n$-dimensional L\'evy process $V$ is a $VG^n(b,\bfmu,\Sigma)$-process if and only if its characteristic exponent has the form (see~\cite{BKMS16}, their Formula~(2.9))
\begin{equation}\label{charexpoVG}
\Psi_V(\bftheta)\,=\,-b\ln\{(b-\rmi\skal{\bfmu}{\bftheta}+(1/2)\|\bftheta\|^2_\Sigma)/b\}\,,\quad \bftheta\in\RR^n\,.\end{equation}
Alternatively, a drift-less $FV^n$-process $X$ is a $VG^n(b,\bfmu,\Sigma)$-process if and only if its L\'evy measure satisfies $\XXX=\VVV_{b,\bfmu,\Sigma}$ for some $b>0,\bfmu\in\RR^n$, a covariance matrix $\Sigma\in\RR^{n\times n}$ and $B\sim BM^n(\bfmu,\Sigma)$, where
\begin{equation}\label{levyVGgen}\VVV_{b,\bfmu,\Sigma}(\rmd \bfy)\,:=\,\int_{(0,\infty)}\PP(B(g)\!\in\! \rmd\bfy)be^{-bg}\,\rmd g/g\,.\end{equation}
This follows from the formula of the L\'evy measure under univariate subordination (see~\cite{s},~his Formula (30.8)).

If, in addition, $\Sigma$ is invertible, then $\VVV_{b,\bfmu,\Sigma}$ is absolutely continuous
with respect to the Lebesgue measure $\rmd\bfv$ on $\RR^n_*$, having L\'evy density $\nu_{b,\bfmu,\Sigma}(\bfv):=(\rmd\VVV_{b,\bfmu,\Sigma}/\rmd\bfv)(\bfv)$, $\bfv\in\RR^n_*$, where (see~\cite{BKMS16}, their Formula~(2.11))
\begin{equation}\label{defgVGbmuSigma}
\nu_{b,\bfmu,\Sigma}(\bfv)\,=\,
\frac{2b\exp\{\skal\bfmu\bfv_{\Sigma^{-1}}\}}{(2\pi)^{n/2}|\Sigma|^{1/2}\|\bfv\|^{n}_{\Sigma^{-1}}}\;
\KKKK_{n/2}\big\{(2b+\|\bfmu\|^2_{\Sigma^{-1}})^{1/2}\|\bfv\|_{\Sigma^{-1}}\,\big\}\,,\end{equation}
$\KKKK_{\rho}(r):=r^{\rho}K_{\rho}(r)$, $\rho\!\ge\!0$, $r\!>\!0$, and $K_\rho$ is a modified Bessel function of the second kind~(see~\cite{BKMS16} and~\cite{GrRy96}, their Equation~(2.12) and their Equation~(3.471)--9, respectively).\\[1mm]
\noindent{\bf Multivariate time parameter.}~If $X,Y\sim L^n$ are independent $n$-dimensional L\'evy processes, then the $2n$-dimensional process $Z=(X,Y)\sim L^{2n}(\bfm,\Sigma,\ZZZ)$ is a L\'evy process in $\RR^{2n}$, for some $\bfm\in\RR^{2n}$, $\Sigma\in\RR^{2n\times 2n}$, and L\'evy measure $\ZZZ$ on $\RR^{2n}_*$. Our notation extends from $\RR^n$ to $\RR^{2n}$ in a canonical way; in particular, $\|\cdot\|$ and $\DD$ may refer to the Euclidean norm and the Euclidean unit ball in $\RR^{n}$ as well as in $\RR^{2n}$, respectively.

As a first step, we evaluate an $n$-dimensional L\'evy process $X=(X_1,\dots, X_n)$, indexed by univariate time $t$, at
multivariate time points ${\bf t}=(t_1,\dots ,t_n)\in[0,\infty)^n$. The result is an infinitely divisible row vector
$X(\bft)\!=\!(X_1(t_1),\dots,X_n(t_n))$. To provide formulae for the associated characteristics, we introduce an operation $\tr$ as an outer product.

For $\bft=(t_1,\dots, t_n)\in[0,\infty)^n$, $\bfmu=(\mu_1,\dots \mu_n)\in \RR^n$ and $\Sigma=(\Sigma_{kl})\in\RR^{n\times n}$, introduce~$\bft\tr\bfmu\in\RR^n$ and~$\bft \tr \Sigma=(\bft\tr\Sigma_{kl})\in\RR^{n\times n}$ by
\begin{equation}\label{bfmux}\bft\tr\bfmu:=(t_1\mu_1,\dots,t_n\mu_n)\,,\quad
(\bft\tr  \Sigma)_{kl}:=\Sigma_{kl}\; (t_k\wedge t_l)\,,\quad 1\le k,l\le n\,.\end{equation}
Choose an ordering $t_{(1)}\le{\dots}\le t_{(n)}$ of the components with associated permutation $\langle(1),\dots,(n)\rangle$ and spacings $\Delta t_{(k)}:=t_{(k)}-t_{(k-1)}$ for
$1\!\le\! k\!\le\! n$, $t_{(0)}:=0$. If $\XXX$ is a L\'evy measure, thus a Borel measure on $\RR^n$
satisfying~\eqref{Piintegrab}, so is $\bft\tr\XXX$, defined by
\begin{equation}\label{VVVsubord}\bft\tr\XXX:=\sum_{k=1}^n \Delta t_{(k)}\XXX_{\{(k),\dots,(n)\}}\,.\end{equation}

We introduce $\bfc(\bft,\XXX)\in\RR^n$ by setting\begin{equation}\label{defbfc}
\bfc\,:=\,\bfc(\bft,\XXX)\,:=\,\sum_{k=2}^n \Delta t_{(k)}\!\!\int_{\DD^C}\!\bfpi_{\{(k),\dots,(n)\}}(\bfx)\,\eins_\DD(\bfpi_{\{(k),\dots,(n)\}}(\bfx))\XXX(\rmd \bfx).
\end{equation}
As~\eqref{Piintegrab} is satisfied for a L\'evy measure $\XXX$, $\bfc(\bft,\XXX)$ is a well-defined $n$-dimensional row vector, and it acts as a compensation term.

We provide formulae for the characteristics of $X(\bft)$ (see Subsection~\ref{subsecproofI} for a proof).
\begin{proposition}\label{propmultiparameter} For $\bft\!=\!(t_1,\dots,t_n)\!\in\![0,\infty)^n$ and $X\!\sim\! L^n(\bfmu,\Sigma,\XXX)$ with $\Psi$ as in~\eqref{0.1},\label{lemmutliple}
	the vector $X(\bft)\!=\!(X_1(t_1),\dots,X_n(t_n)))\in\RR^n$ is infinitely divisible
	with $\Phi_{X(\bft)}(\bftheta)\!=\!\EE\exp(\rmi\skal\bftheta {X(\bft)})\!=\!\exp((\bft\tr\Psi)(\bftheta))$, $\bftheta\!\in\!\RR^n$, where
\begin{eqnarray}
\lefteqn{\hspace*{-3em}(\bft\tr \Psi)(\bftheta)\,:=\,\sum_{k=1}^n\Delta t_{(k)} \Psi(\bfpi_{\{(k),\dots,(n)\}}(\bftheta))}&&\label{multiexponent1}\\
&&\hspace{-1cm}\,=\,\rmi\skal{\bft\tr\bfmu+ \bfc}\bftheta-\frac 12\;\|\bftheta\|^2_{\bft\tr\Sigma}+\int_{\RR_*^n}\big(e^{\rmi\skal\bftheta \bfx}-1-\rmi\skal\bftheta \bfx \eins_{\DD}(\bfx)\big)\,\bft\tr\XXX(\rmd \bfx)\,.\label{multiexponent2}
\end{eqnarray}
\end{proposition}
\begin{remark}\label{rem1counterexamplBBT2T}
	If $B$ is a standard Brownian motion and $I$ is the identity function, then
	$(I,2I)$ is a subordinator and $(B,B)$ is a L\'evy process, but $(B,B)\circ(I,2I)$ is not a L\'evy process, as follows from Proposition~\ref{propindepindist}(iii) below. Though it is a Gaussian
	process, $(B,B)\circ(I,2I)$ is not a Brownian motion.\halmos \end{remark}
\noindent{\bf Weak subordination.}~If $T\sim S^n(\bfd,\TTT)$ is a subordinator, then
we may write $T=I\bfd+S$, where $S\sim S^n({\bf 0},\TTT)$ is a pure jump subordinator and $I\bfd$ is a deterministic subordinator. Suppose $X\sim L^n(\bfmu,\Sigma,\XXX)$ is the candidate for a subordinate. If $Y$ is another L\'evy process with
$Y\sim L^n(\bfd\tr\bfmu\!+\!\bfc,\bfd\tr\Sigma,\bfd\tr\XXX)$,
we get from Proposition~\ref{propmultiparameter} that $X(t\bfd)\eqd Y(t)$ for all fixed $t\!\ge\! 0$. Any other L\'evy process with this property must have the same characteristics, and in the case of deterministic subordination
the law of weak subordination is determined.

On the other hand, if $\bfd\!=\!{\bf 0}$ and $T\!=\!S$, we may perceive the subordinated process
as a $2n$-dimensional L\'evy process $Z\!=\!(S,Y)$ in time-space $[0,\infty)\times\RR^n_*$, and the jumps
of $Y$ should have conditional laws $(\Delta Y|\Delta T=\bft) \sim \PP(X(\bft)\in\cdot)$. This notion is consistent with traditional subordination~(see~\cite{BPS01} and \cite{s}, their Formula (3.12) and his Formula (30.8), respectively) as illustrated by \eqref{levyVGgen} in the context of $VG^n$-processes. Equivalently, the jumps of the joint process $Z$ form a marked point process,
with marks in time-space $[0,\infty)\times\RR^n_*$ determined, conditionally independently, based on the
points of a Poisson point process in time-time $[0,\infty)\times[0,\infty)^n_*$, with intensity measure $\rmd t\otimes \TTT$. Summing up those jumps along $t\ge 0$, possibly with a compensation term, generates a pure-jump L\'evy process $Z$ with values in time-space $[0,\infty)^n_*\times \RR^n_*$. Using pure-jump subordinators, the law of weak subordination is thus determined in time-space.

Traditional subordination is consistent with the superposition of independent subordinators such as $T=I\bfd+S$, and the law of strongly subordinated processes, when defined as L\'evy processes,
is determined by convolution (see~\cite{BKMS16}, their Proposition~4.1). Together with Proposition~\ref{propsupervisexpo}, this imposes a final and determining constraint on the law of weak subordination.

We are now prepared to introduce subordination in the weak and semi-strong senses.
\begin{definition}\label{defsupervision} Let $X\sim L^n(\bfmu,\Sigma,\XXX)$ and $T\sim S^n(\bfd,\TTT)$. A process $Z$ is called {\em a subordinator $T$ subordinating $X$ in the weak sense},
meaning that $Z\eqd(T,X\su T)$, whenever
$Z=(Z_1,Z_2)\sim L^{2n}(\bfm,\allowbreak\Theta,\ZZZ)$ is a L\'evy process
with the characteristics determined by $\bfm\!=\!(\bfm_1,\bfm_2)$,
$\bfm_1,\bfm_2\in\RR^{n}$,
\begin{align}
\bfm_1&=\bfd+\int_{[0,\infty)_*^n} \bft\,\PP\big((\bft,X(\bft))\in\DD\big)\;\TTT(\rmd \bft)\,,\label{defmu1}\\
\bfm_2&=\bfc(\bfd,\XXX)+\bfd\tr\bfmu+\int_{[0,\infty)^n_*}\EE[X(\bft)\;\eins_{\DD}(\bft,X(\bft))]\,\TTT(\rmd \bft)\,,\label{defmu2}\\
\Theta&=\left(\begin{array}{cc}\bfnull&\bfnull\\\bfnull&{\bfd\tr\Sigma}\end{array}\right)\,,\label{SigmaY}\\
\ZZZ(\rmd \bft,\rmd \bfx)&=(\bfdelta_{{\bf 0}}\!\otimes\! (\bfd\!\tr\!\XXX))(\rmd\bft,\rmd\bfx)\!+\!\eins_{[0,\infty)^n_*\times\RR^n}\PP(X(\bft)\!\in\!\rmd\bfx)\TTT(\rmd \bft)\,.\label{PliY}
\end{align}
We call $Z=(Z_1,Z_2)$ {\em a subordinator $T$ subordinating~$X$ in the semi-strong sense}, whenever, simultaneously, $Z_1\!=\!T$ are indistinguishable and $Z\!\eqd\!(T,X\odot~T)$.\halmos
\end{definition}
Such a process $Z$ exists and is a L\'evy process whenever the specifying characteristics are as in~\eqref{defmu1}--\eqref{PliY}. The main difficulty is to show that
$\PP(X(\bft)\in\rmd\bfx)\;\TTT(\rmd \bft)$ is a L\'evy measure. Semi-strong subordination is then always possible on augmented probability spaces, and it relies on marking the point process associated to the jumps of $T$ (see Subsection~\ref{subsecproofII} for a proof).
\begin{theorem}\label{thmexist} Let $X$ and $T$ be as in Definition~\ref{defsupervision}.\\
(i)~There exists a L\'evy process $Z=(Z_1,Z_2)\sim L^{2n}(\bfm,\Theta,\ZZZ)$ with $\bfm,\Theta,\ZZZ$ as specified in \eqref{defmu1}--\eqref{PliY}.\\
(ii)~On an augmentation of the probability space carrying $T$, there exists an $n$-dimensional
L\'evy process $Z_2$ such that $(T,Z_2)$ is a subordinator $T$ subordinating $X$ in the semi-strong sense.\\
(iii)~If both $\bfd\!=\!{\bf 0}$ and $\int_{[0,1]_*^n}\|\bft\|^{1/2}\;\TTT(\rmd\bft)\!<\!\infty$ hold, then $Z$ is a drift-less $FV^{2n}$-process.
\end{theorem}
\begin{remark}\label{rem2counterexamplBBT2T}
	Let $B,B^*$ be independent standard Brownian motions, and $I$ be the identity function. $(I,2I)$ is a subordinator, and $(B,B)$ is a L\'evy process, but $(B,B)\circ(I,2I)$ is not (see Remark~\ref{rem1counterexamplBBT2T}).
As easily verified from~\eqref{defmu1}--\eqref{PliY}, the process
	$(I,2I,B,B+B^*)$ is the subordinator $(I,2I)$ subordinating $(B,B)$ in the semi-strong sense. \halmos\end{remark}
\begin{remark}Recall $\bfe=(1,\dots,1)\in\RR^n$. Let $G\sim\Gamma_S(b)=S^1(0,\GGG_b)$.
In the traditional notion of $VG^n$, $G$ is the sole time change for all the components of an $n$-dimensional Brownian motion. Consistent with our notion of multivariate subordination, we replace $G$ with $G\bfe$. Note $G\bfe=(G,\dots,G)$ is an $n$-dimensional drift-less subordinator with indistinguishable
components and $G\bfe\sim S^n({\bf 0},\GGG_b\circ (I\bfe)^{-1})$.
\halmos\end{remark}
\begin{remark}
The $VG^n$-process in~\cite{MaSe90} uses $n$-dimen\-sional Brownian motion as its subordinate and a univariate standard gamma process as its subordinator.
The $VG^n$-model gives a restrictive dependence structure, where components cannot have idiosyncratic time changes and must have equal kurtosis. These last two deficiencies have been addressed by Luciano and Semeraro's~\cite{LS10,Se08} variance-$\bfalpha$-gamma ($V\bfalpha G$) process by the use of a $\bfalpha$-gamma subordinator.\halmos
\end{remark}
\noindent{\bf $\bfalpha$-gamma subordinator.}~Assume $n\!\ge\!2$. Let $a,b\!>\!0$, $\bfalpha\!=\!(\alpha_1,\dots,\alpha_n)\!\in\!(0,\infty)^n$
such that $b\!>\!a\alpha_k$ for $k=1,\dots,n$. Introduce $\beta_k:=(b\!-\!a\alpha_k)/\alpha_k$, and let $G_0,\dots,G_{n}$ be independent gamma subordinators such that $G_{0}\sim\Gamma_S(a,b)$, $G_k\sim\Gamma_S(\beta_k,b/\alpha_k)$, $1\!\le\! k\!\le\! n$.

We refer to an $n$-dimensional subordinator $T\sim \bfalpha G^n(a,b,\bfalpha)$, as an
{\em$\bfalpha$-gamma ($\bfalpha G$) subordinator}~\cite{Se08}, provided
\begin{equation}\label{semeraro2}
T=(T_1,\dots,T_n)\,\eqd G_0 \bfalpha+(G_1,\dots,G_n)\,.
\end{equation}
As perceived in~\cite{Se08}, the components of $T$ are univariate {\em standard gamma} subordinators $
T_k\sim\Gamma_S(b/\alpha_k)$, $1\!\le\! k\!\le\! n$. Further, an $n$-dimensional drift-less subordinator $T$ with L\'evy measure $\TTT$
is an $\bfalpha$-gamma subordinator with parameters $a,b,\bfalpha$ if and only if, with $\beta_1,\dots,\beta_n$ as above, \begin{align}\label{levyalphasubord}\TTT=
\int_{(0,\infty)} \bfdelta_{g\bfalpha}\;\GGG_{a,b}(\rmd g)+
\sum_{k=1}^{n}\bfdelta_0^{\otimes (k\!-\!1)}\otimes\GGG_{\beta_k,b/\alpha_k}\otimes \bfdelta_0^{\otimes(n\!-\!k)}\,.\end{align}
\noindent{\bf Strong variance-$\bfalpha$-gamma processes.}~Assume $n\!\ge\!2$. 
Let $\bfmu=(\mu_1,\dots,\mu_n)\!\in\!\RR^n$ and $\Sigma=\diag(\Sigma_{11},\dots,\Sigma_{nn})\!\in\![0,\infty)^{n\times n}$ be a {\em diagonal} matrix.

For independent $B$ and $T$, if $Y\eqd B\circ T$, where $B\sim BM^n(\bfmu,\allowbreak \Sigma)$ is a Brownian motion with independent components and $T\sim \bfalpha G^n(a,b,\bfalpha)$ is an $\bfalpha G$-subordinator,
then we call $Y$ a {\em $($strong$)$ variance-$\bfalpha$-gamma process} with parameters $a,b,\bfalpha,\bfmu,\Sigma$, in short,
\begin{equation}\label{defstrongValphaG} Y\eqd B\circ T\sim V\bfalpha G^n(a,b,\bfalpha,\bfmu,\Sigma)=BM^n(\bfmu,\Sigma)\circ \bfalpha G_S^n(a,b,\bfalpha)\,.
\end{equation}
\begin{remark}The Brownian motion subordinate must have {\em independent} com\-ponents, which re\-stricts the de\-pendence struc\-ture.
In our {\em weak form\-ulation} of the $V\bfalpha G$-process, the subordinate is a Brownian motion with possibly {\em correlated} components. Our $WV\bfalpha G$-process has a wider range of de\-pendence struc\-tures, while being pa\-rsi\-moniously para\-metrised, each component has both common and idio\-syncratic time changes, it has $VG$-marginals with in\-dependent levels of kur\-tosis, with the jump measure having full support.~\halmos
\end{remark}
\noindent{\bf Weak variance-$\bfalpha$-gamma processes.}~Assume $n\!\ge\!2$. 
Let $\bfmu=(\mu_1,\dots,\mu_n)\!\in\!\RR^n$ and $\Sigma=(\Sigma_{kl})\!\in\!\RR^{n\times n}$ be an {\em arbitrary} covariance matrix.

Whenever $Y\eqd B\odot T$, where $B\sim BM^n(\bfmu,\Sigma)$ is Brownian motion, and $T\sim \bfalpha G^n(a,b,\bfalpha)$ is an $\bfalpha G$-subordinator,
then we call $Y$ a {\em weak variance-$\bfalpha$-gamma process} with parameters $a,b,\bfalpha,\bfmu,\Sigma$, in short,
\begin{equation}\label{defweakValphaG} Y\eqd B\odot T\sim WV\bfalpha G^n(a,b,\bfalpha,\bfmu,\Sigma)=BM^n(\bfmu,\Sigma)\odot\bfalpha G_S^n(a,b,\bfalpha)\,.
\end{equation}
We derive the joint L\'evy measure $\ZZZ$ of the pair $Z=(T,B\odot T)$. Let $N(\rmd \bfx|\bfmu,\Sigma)$ be the normal law with mean $\bfmu$ and covariance matrix $\Sigma$. As
$B(\bft)\sim N(\rmd\bfx|\allowbreak\bft\tr\bfmu,\bft\tr\Sigma)$ for fixed $\bft\in[0,\infty)_*^n$, by Proposition~\ref{propmultiparameter},
it follows from~\eqref{PliY} and~\eqref{levyalphasubord} that for a Borel set $A\subseteq \RR^{2n}_*$,
\begin{eqnarray}
\lefteqn{\hspace*{-1em}\ZZZ(A)\,=\,\int_{(0,\infty)\times\RR^n}\eins_A(g\bfalpha,\bfx) N(\rmd \bfx|g \bfalpha \tr \bfmu,g\bfalpha\tr\Sigma)\GGG_{a,b}(\rmd g)}&&\nonumber\\
&&\quad+\sum_{k=1}^n\int_{(0,\infty)\times\RR}\eins_{A}(g\bfe_k,x_k\bfe_k)\,N(\rmd x_k|g\mu_k,g\Sigma_{kk})\;\GGG_{\beta_k,b/\alpha_k}(\rmd g)\,.\label{levyalphaVII}
\end{eqnarray}
Formula~\eqref{levyalphaVII} tells us that
$T$ and $B\odot T$ jump together. As a result, weakly subordinated Brownian motion resembles the jump behaviour of a subordinated Brownian motion. Like strong $V\bfalpha G$-processes in~\cite{Se08}, $WV\bfalpha G$-processes jump in two different ways: either the components jump independently
of each other together with one of the subordinators $G_1,\dots,G_n$, or
the components jump together with the subordinator~$G_0$.

If, in addition, $\Sigma$ is a diagonal matrix, and~\eqref{levyalphaVII} is projected on space we recover
the formulae of the strong variance-$\bfalpha$-gamma process $B\circ T$, as derived in~\cite{LS10} (see their Theorem~1.1).
\section{Properties of Weak Subordination} \label{secprop}
Let $T\sim S^n(\bfd,\TTT)$ and $X\sim L^n(\bfmu,\Sigma,\XXX)$ be candidates for a subordinator and
subordinate in the weak or semi-strong subordination of Definition~\ref{defsupervision}.

We provide a formula for the characteristic exponent.
\begin{proposition}\label{propcharexpo}
	$Z\eqd (T,X\odot T)$ holds in the weak sense if and only if
	for all $\bftheta=(\bftheta_1,\bftheta_2)$ with $\bftheta_1,\bftheta_2\in\RR^n$, the
	characteristic exponent of $Z$ is
	\begin{equation}\label{propsupervisexpo}
\Psi_Z(\bftheta)\!=\!\rmi\skal{\bfd}{\bftheta_1}+(\bfd\tr\Psi_{X})(\bftheta_2)+\int_{[0,\infty)^n_*}\!(\Phi_{(\bft,X(\bft))}(\bftheta)-1)\TTT(\rmd \bft)\,.
	\end{equation}
	Here, $(\bfd\tr\Psi_{X})(\bftheta_2)$ is defined as in~\eqref{multiexponent2}, but with $(\bftheta,\bft)$ replaced by $(\bftheta_2,\bfd)$.
\end{proposition}
\noindent{\bf Proof.}~Let $\bftheta=(\bftheta_1,\bftheta_2)$ with $\bftheta_1,\bftheta_2\in\RR^n$. Combining~\eqref{PiFVintegrab} and Lemma~\ref{lemcharlinbounf} below yields $\TTT$-integrability of $\bft\mapsto \EE\exp\{\rmi\skal{\bftheta}{(\bft,X(\bft))}\}-1$ with
\begin{align}\label{idcharexpo}
&\hspace{-1cm}\int_{[0,\infty)_*^n\times\RR^n}\big(e^{\rmi\skal\bftheta {(\bft,\bfx)}}\!-\!1\!-\!\rmi\skal\bftheta{(\bft,\bfx)} \eins_\DD(\bft,\bfx)\big)\,\PP(X(\bft)\in\rmd \bfx)\,\TTT(\rmd \bft)\\
=&{}-\rmi\int_{[0,\infty)^n_*}(\EE[\skal{\bftheta_2}{X(\bft)}\eins_{\DD}(\bft,X(\bft))]+\skal{\bftheta_1}\bft \PP((\bft,X(\bft))\in\DD)\;\TTT(\rmd \bft)\nonumber\\
&{}+\int_{[0,\infty)^n_*}\!\big(\Phi_{(\bft,X(\bft))}(\bftheta)-1\big)\TTT(\rmd \bft)\nonumber
\end{align}
Plainly, $Z\eqd (T,X\odot T)$ if and only if $Z$ has characteristic triplet~\eqref{defmu1}--\eqref{PliY}, and this is true if and only if~\eqref{propsupervisexpo} does, as follows from~\eqref{idcharexpo}.\halmos\\

Apart from determining the distribution of the time-and-space projected processes, the next proposition states that, in analogy with traditional subordination~\cite{BPS01,s},
weak subordination is consistent with projections.
\begin{proposition}\label{propprops} If $Z=(Z_1,Z_2)$ is a L\'evy process with $n$-dimensional
	components $Z_1$ and $Z_2$ such that $Z=(Z_1,Z_2)\eqd (T,X\odot T)$ holds in the weak sense, then we must have~$Z_1\eqd T$ as well as $Z_2\sim L^n(\bfm_{2},\Theta_2,\ZZZ_2)$ with
	\begin{align}
	\bfm_2&=\bfc(\bfd,\XXX)\!+\!\bfd\!\tr\!\bfmu
	\!+\!\int_{[0,\infty)^n_*}\EE[X(\bft)\eins_{\DD}(X(\bft))]\;\TTT(\rmd \bft),\label{muY2}\\
\Theta_2&=\bfd\tr\Sigma\label{SigmaY2}\\
	\ZZZ_2(\rmd\bfx)&=\bfd\tr\XXX(\rmd\bfx)+\,\int_{[0,\infty)^n_*}\PP(X(\bft)\in\rmd\bfx)\;\TTT(\rmd \bft)\label{PliY2}\,.
	\end{align}
In addition, if $J\subseteq \{1,\dots,n\}$, we have $(T\bfpi_J,(X\odot T)\bfpi_J)\eqd (T\bfpi_J,(X\bfpi_J)\odot (T\bfpi_J))$ and, particularly, $(T_k,(X\odot T)_k)\eqd (T_k,X_k\odot T_k)$ for $1\!\le\! k\!\le\! n$.
\end{proposition}
\noindent{\bf Proof.}\label{subsecproofpropprops}~Let $Z=(Z_1,Z_2)\sim L^{2n}(\bfm,\Theta,\ZZZ)$ with $\bfm=(\bfm_1,\bfm_2),\Theta,\ZZZ$ as specified in \eqref{defmu1}--\eqref{PliY}. For $\bftheta_1,\bftheta_2\in\RR^n$ it is straightforwardly checked that $\Psi_{Z}(\bftheta_1,{\bf 0})=\Psi_T(\bftheta_1)$ and $\Psi_{Z}({\bf 0},\bftheta_2)=\Psi_{Z_2}(\bftheta_2)$, giving $Z_1\eqd T$ and $Z_2\sim L^n$ with characteristics matching those in~\eqref{muY2}--\eqref{PliY2}.

Without loss of generality, assume $J\neq\emptyset$ and
set $\bfpi:=\bfpi_J$. It suffices to show that, for all $\bftheta_1,\bftheta_2\in\RR^n$,
\begin{equation}\label{eqproj}
\Psi_{(T\bfpi,(X\odot T)\bfpi)}(\bftheta_1,\bftheta_2)=\Psi_{(T\bfpi,(X \bfpi)\odot (T\bfpi))}(\bftheta_1,\bftheta_2)\,.
\end{equation}
By noting $\Psi_{(T\bfpi,(X\odot T)\bfpi)}(\bftheta_1,\bftheta_2)=\Psi_{(T,X\odot T)}(\bftheta_1\bfpi,\bftheta_2\bfpi)$, the LHS in~\eqref{eqproj} matc\-hes the RHS in~\eqref{propsupervisexpo}, but with $(\bftheta_1,\bftheta_2)$ replaced with $(\bftheta_1\bfpi,\bftheta_2\bfpi)$.
The RHS in~\eqref{eqproj} equals the RHS in~\eqref{propsupervisexpo} with $(T,X)$ replaced with $(T\bfpi,X\bfpi)$. To prove the identity in~\eqref{eqproj}, it thus suffices to compare the three terms
occurring on both sides in \eqref{eqproj}, respectively.

The projected process $T\bfpi$ is an $n$-dimensional subordinator with drift $\bfd\bfpi$ and L\'evy measure $\TTT\circ \bfpi^{-1}$. Consequently, the first term on both sides in~\eqref{eqproj} are equal as $\langle\bftheta_1\bfpi,\bfd\rangle=\langle\bftheta_1,\bfd\bfpi\rangle$. The second identity,
$\bfd\tr\Psi_X(\bftheta_2\bfpi)=(\bfd\bfpi)\tr\Psi_{X\bfpi}(\bftheta_2)$,
follows from Proposition~\ref{propmultiparameter} as $\langle\bftheta_2\bfpi,X(\bfd)\rangle=\langle\bftheta_2,(X(\bfd))\bfpi\rangle=
\langle\bftheta_2,(X\bfpi)(\bfd\bfpi)\rangle$. The third identity follows from the transformation theorem by recalling that $\TTT\circ \bfpi^{-1}$ is the L\'evy measure of the projected process $T\bfpi$,
and by $(X \bfpi)(\bft)=(X\bfpi)(\bft\bfpi)$ and $\langle(\bftheta_1\bfpi,\bftheta_2\bfpi),(\bft,X(\bft))\rangle=\skal {(\bftheta_1,\bftheta_2)}{(\bft\bfpi,(X\bfpi)(\bft\bfpi)}$, $\bft\in[0,\infty)^n$, as they imply the crucial identity
\[\int_{[0,\infty)_*^n} \big(\Phi_{(\bft,X(\bft))}(\bftheta_1\bfpi,\bftheta_2\bfpi)-1\big)\TTT(\rmd\bft)
=\int_{[0,\infty)_*^n}\big(\Phi_{(\bft,(X\bfpi)(\bft))}(\bftheta)-1\big)\TTT\circ \bfpi^{-1}(\rmd\bft)\,.\]
\halmos
\begin{remark}\label{remmargconsVaG} Weak subordination is consistent with projections to coordinates by Proposition~\ref{propprops}. Suppose $Y=(Y_1,\dots,Y_n)\eqd B \odot T\sim WV\bfalpha G^n(a,b,\bfalpha,\bfmu,\Sigma)$ in~\eqref{defweakValphaG}, where
$B=(B_1,\dots,B_n)\sim BM^n(\bfmu,\Sigma)$ and $T=(T_1,\dots,T_n)\sim\bfalpha G^n(a,b,\bfalpha)$. Assume that $B$ and $ T$ are independent, then $Y$ has $VG^1$-components. Thus, $Y$ has the same marginal distributions as a strong $V\bfalpha G^n(a,b,\bfalpha,\bfmu,\Sigma)$-process~\cite{Se08} because\begin{equation}\label{VGcomps}
Y_k\eqd(B\odot T)_k\eqd B_k\odot T_k\eqd B_k\circ T_k \sim VG^1(b/\alpha_k,\mu_k,\Sigma_{kk})\,,\quad 1\!\le\!k\!\le\!n\,.\end{equation}
\halmos\end{remark}
Weak and semi-strong subordination extends traditional subordination.
\begin{proposition}\label{propsupextendssub} Let $T,X$ be independent. If either $T$ has indistinguishable components or $X$ has independent components, then $(T,X\circ T)\eqd(T,X\odot T)$ in the semi-strong sense.
\end{proposition}
\noindent{\bf Proof.}\label{subsecproofpropsupextendssub}~We extend $\skal {\bf z}{\bf w}:=\sum_{k=1}^n z_kw_k$ to ${\bf z},{\bf w}\in\CC^n$. We avoid conjugation.

As we assumed $T$ and $X$ to be {\em independent} processes, we get from Proposition~\ref{propmultiparameter} by conditioning on $T$ that, for
$\bftheta=(\bftheta_1,\bftheta_2),\bftheta_1,\bftheta_2\in\RR^n$,
\begin{equation}\label{charsubord}\Phi_{(T(1),X(T(1)))}(\bftheta)\,=\,\EE\exp\{\rmi \skal{\bftheta_1}{T(1)}+(T(1)\tr\Psi_{X})(\bftheta_2)\}\,.
\end{equation}
{\it Univariate subordination.}~$T,X$ are independent with $T=R\bfe$ with $R\sim S^1(d,\RRR)$ and $\bfe=(1,\dots,1)\in\RR^n$.
We have $\bfc={\bf 0}$ in \eqref{defbfc}. Note $\skal{\bftheta_1}{T(1)}=R(1)\skal{\bftheta_1}{\bfe}$ and $(T(1)\tr\Psi_{X})(\bftheta_2)=R(1)(\bfe\tr\Psi_{X}(\bftheta_2))$
in~\eqref{charsubord}. Noting $\Re z\ge 0$ for
$z:=-\rmi\skal{\bftheta_1}{\bfe}-\bfe\tr\Psi_{X}(\bftheta_2)$, we get from~\eqref{charsubord} that
$\Psi_{(T,X\circ T)}(\bftheta)=-\Lambda_R(z)$, where $\Lambda_R(z):=d z+\int_{(0,\infty)}(1-e^{-zr})\;\RRR(\rmd r)$.
The RHS matches~\eqref{propsupervisexpo}, and $T$ subordinates $X$ in the semi-strong sense.\\
{\it Multivariate subordination.} Let $T,X_1,X_2,\dots,X_n$ be independent. Particularly,~$\Sigma$ is a diagonal matrix and $\XXX=\sum_{k=1}^n\XXX_{\{k\}}$. If $\bft=(t_1,\dots,t_n)\in[0,\infty)^n$, $J_{(m)}:=\{(m),\dots,(n)\}$, $1\!\le\! m\!\le\! n$, then~\eqref{VVVsubord} becomes $\bft\tr\XXX=\sum_{k=1}^n t_k \XXX_{\{k\}}$ as
\[\sum_{m=1}^n\Delta t_{(m)}\Big\{\sum_{k=1}^n\XXX_{\{k\}}\Big\}_{J_{(m)}}
	=\sum_{k=1}^n\Big\{\sum_{m=1}^n\Delta t_{(m)} \eins_{J_{(m)}}(k)\Big\}\XXX_{\{k\}}\,.\]
Note $\bfc={\bf 0}$ in \eqref{defbfc} because for $\emptyset\neq J\subseteq\{1,\dots,n\}$, $1\!\le\! k\!\le\! n$,
\[\int_{\DD^C} \bfx\bfpi_{J}\,\eins_{\DD}(\bfx\bfpi_{J})\,\XXX_{\{k\}}(\rmd\bfx)
=\eins_{J}(k)\,\bfe_k\int_{\DD^C} x\,\eins_{\DD}(x)\XXX_{\{k\}}(\rmd x)= {\bf 0}\,.\]
Recalling the diagonal form of $\Sigma$ yields $\skal{\bftheta_2(\bftheta_2\tr\Sigma)}{\bft}=\|\bftheta_2\|^2_{\bft\tr\Sigma}$ for $\bft\in[0,\infty)^n,\bftheta_2\in\RR^n$. Also, $\Re{\bf z}\in[0,\infty)^n$, for $\bftheta_1,\bftheta_2\in\RR^n$ and
\[{\bf z}:=\frac 12 \bftheta_2(\bftheta_2\tr\Sigma)
-\rmi(\bftheta_1+\bftheta_2 \tr\bfmu)-\sum_{k=1}^n\!\bfe_k\!\!\int_{\RR^n_*}(e^{\rmi\skal{\bftheta_2}\bfx}-1-\rmi\skal{\bftheta_2}\bfx\!\eins_{\DD}(\bfx))\XXX_{\{k\}}(\rmd\bfx)\,.\]
By~\eqref{charsubord}, note $\Psi_{(T,X\circ T)}(\bftheta)=-\Lambda_T(\bfz)$, where $\Lambda_T(\bfz):=\skal{\bfz}{\bfd}+\int_{[0,\infty)_*^n}(1-e^{-\skal {\bfz}\bft})\;\TTT(\rmd \bft)$.
As RHS matches~\eqref{propsupervisexpo}, $T$ subordinates $X$ in the semi-strong sense.\halmos
\begin{remark}\label{remVaGweakstrong}
Suppose $Y\eqd B\odot  T\sim WV\bfalpha G^n(a,b,\bfalpha,\bfmu,\Sigma)$ in~\eqref{defweakValphaG}.
If $\Sigma$  is of diagonal form, then $B$ is a Brownian motion with independent increments.
Assume $B$ and $T$ are independent. Proposition~\ref{propsupextendssub} states that $Y\eqd B\odot T\eqd B\circ T\sim V\bfalpha G^n(a,b,\bfalpha,\bfmu,\Sigma)$ in~\eqref{defstrongValphaG}.
Within the general class of $n$-dimensional L\'evy processes, the $WV\bfalpha G$-class is thus a proper extension of the strong $V\bfalpha G$-class.\halmos
\end{remark}
\noindent{\bf Monotone case.}~If its standard assumptions are violated, then
traditional subordination may fail to create L\'evy processes. Curiously, weak subordination overcomes this problem in the {\em monotone} case.
\begin{proposition}\label{propmonocase} Suppose $T,X$ be independent while $Z\eqd(T,X\odot T)$ in the weak sense.
If $T$ has monotone components $T_1\le{\dots}\le T_n$, then $Z(t)\eqd (T(t),X(T(t)))$ for all fixed $t\ge 0$.
\end{proposition}
\noindent{\bf Proof.}\label{subsecproofpropmonocase} Set $[0,\infty)^n_\le:=\{\bft=(t_1,\dots,t_n)\in[0,\infty)^n:\,t_1\!\le\! {\dots}\!\le\! t_n\}$.
For $1\!\le\! k\!\le\! n$, let $\Sigma_{k}=(\Sigma_{k,ij})\in\RR^{n\times n}$ be defined by $\Sigma_{k,ij}:=\Sigma_{ij}\eins\{i \wedge j\ge k\}$ for $1\!\le\! i,j\!\le\! n$. For $\bft=(t_1,\dots,t_n)\in[0,\infty)^n_{\le}$, $1\le k\le n$, let $\Delta t_k:=t_{k}\!-\!t_{k-1}$, $t_0:=0$. The quantities in~\eqref{bfmux}--\eqref{defbfc} are
$\bft\tr\Sigma=\sum_{k=1}^n\Delta t_{k}\Sigma_k$, $\bft\tr\XXX=\sum_{k=1}^n\Delta t_{k}\XXX_{\{k,\dots,n\}}$ and $\bfc=\sum_{k=2}^n\Delta t_{k}\int_{\DD^C} \bfx \bfpi_{\{k,\dots,n\}}\eins_{\DD}(\bfx\bfpi_{\{k,\dots,n\}})\,\XXX(\rmd\bfx)$.

Introduce linear bijections $A,D:\RR^{n}\to\RR^n$ by setting
\[\hspace*{-2cm}\bfx A\,:=\,(x_1,x_1+x_2,x_1+x_2+x_3,\dots,x_1+x_2+\dots +x_n),\]
\[\bfx D\,:=\,(x_1,x_2-x_1,x_3-x_2,\dots,x_n-x_{n-1})\,,\quad\bfx=(x_1,\dots,x_n)\in\RR^n\,.\]
As we assumed $T_1\le{\dots}\le T_n$, $TD\sim S^n(\bfd D,\TTT\circ D^{-1})$ is a subordinator~(see~\cite{s}, his Theorem~24.11).

Let $\bftheta=(\bftheta_1,\bftheta_2)$, $\bftheta_1,\bftheta_2\in\RR^n$. Observe that~$\Re{\bf z}\in[0,\infty)^n$, where
\begin{eqnarray*}
	{\bf z}&:=&-\rmi \bftheta_1A'-\rmi (\bftheta_2\tr\bfmu)A'+\frac 12\sum_{k=1}^n\|\bftheta_2\|^2_{\Sigma_k}\bfe_k\\
	&&\;-\sum_{k=1}^n\int_{\RR^n_*}(e^{\rmi\skal{\bftheta_2}\bfx}-1-\rmi\skal{\bftheta_2}\bfx\!\eins_{\DD}(\bfx))\XXX_{\{k,\dots,n\}}(\rmd\bfx)\,\bfe_k\\
	&&\;-\rmi \sum_{k=2}^n\int_{\DD^C} \skal{\bftheta_2}{\bfx\bfpi_{\{k,\dots,n\}}}\eins_{\DD}(\bfx\bfpi_{\{k,\dots,n\}})\,\XXX(\rmd\bfx)\,\bfe_k\,.
\end{eqnarray*}
As $A=D^{-1}$, note~$\rmi\skal{\bftheta_1}{\bft}+\bft\tr\Psi_{X}(\bftheta_2)=-\skal{{\bf z}}{\bft D}$. Then using the assumption that $T$ and $X$ are independent, and the facts $\bfd\in[0,\infty)^n_\le$ and $\TTT([0,\infty)^n_*\backslash[0,\infty)^n_\le)=0$, it follows that $\Phi_{(T(t),X\circ T(t))}(\bftheta)=\exp\{-t\Lambda_{TD}(\bfz)\}$ for $t\geq 0$, where $\Lambda_{TD}(\bfz):=\skal{{\bf z}}{\bfd D}+\int_{[0,\infty)^n_*}(1-e^{-\skal{{\bf z}}{{\bft}D}})\;\TTT(\rmd \bft)$ matches~\eqref{propsupervisexpo}.\halmos
\begin{remark}\label{rem3counterexamplBBT2T}Let $B,B^*,I$ be the processes specified in Remark~\ref{rem2counterexamplBBT2T} so that $Z\!:=\!((I,2I),(B,B)\odot(I,2I))\eqd(I,2I,B,B+B^*)$. The deterministic subordinator $(I,2I)$ satisfies $I\le 2I$. Proposition~\ref{propmonocase} matches $Z(t)\!\eqd\!(t,2t,B(t),B(2t))$ for all fixed $t\!\ge\!0$.\halmos\end{remark}
\begin{remark}\label{remcounterexamplmonotone}
Suppose $B$ is a standard Brownian motion and $N$ is a Poisson process with unit rate, independent of $B$. Note $\EE[B(t)B(N(t))]\!=\!\EE[t\wedge N(t)]
\!=\!t (1\!-\!e^{-t})$ for $0\!\le\!t\!\le\!1$, which is a nonlinear function in $t$. As a result, $(B,B)\circ (I,N)$ cannot be a L\'evy process, and there
	is no L\'evy process matching $(B,B)\circ (I,N)$ in law in all fixed time points $t\ge 0$. Neither $I\le N$ nor $N\le I$ holds for the subordinator $(I,N)$.
It is verified from~\eqref{defmu1}--\eqref{SigmaY} that
	$((I,N),(B,B)\odot (I,N))\eqd((I,N),(B^*,B\circ N))$ in the semi-strong sense, where $B^*\eqd B$ is independent of $B,N$.\halmos
\end{remark}
\noindent{\bf Ray-subordination.}~Recall $\bfe\!=\!(1,\dots,1)\!\in\!\RR^{n}$, and let $(\bfe,\bfe)=(1,\dots,1)\!\in\!\RR^{2n}$. If $\bfalpha\!\in\![0,\infty)^n$
is a deterministic vector and $R$ is a univariate subordinator, then $T\!:=\!R\bfalpha$
defines an $n$-dimensional subordinator travelling along the deterministic ray $\{r\bfalpha:r\ge 0\}$. We refer to this kind of subordination as ray-subordination. A special case
is strong univariate subordination where the corresponding ray is given by $\{r\bfe:r\ge 0\}$.

Curiously, it is possible
to perceive weak subordination along deterministic rays as univariate subordination of augmented processes.
\begin{proposition}\label{propweakalongray} Let $\bfalpha\in[0,\infty)^n$ be a deterministic vector and $R$ a univariate subordinator.
If $Y$ is a L\'evy process with characteristic exponent $\Psi_Y=\bfalpha\tr \Psi_X$, as in \eqref{multiexponent2}, but with $\bft$ replaced by $\bfalpha$, then we have $(R\bfalpha,X\odot (R\bfalpha))\eqd (I\bfalpha,Y)\odot (R(\bfe,\bfe))$.

If, in addition, $R$ and $Y$ are independent, then $(R\bfalpha,X\odot(R\bfalpha))\eqd(I\bfalpha,Y)\circ(R(\bfe,\bfe))$.
\end{proposition}
\noindent{\bf Proof.}\label{subsecproofpropweakalongray}~
Let $\bftheta\!=\!(\bftheta_1,\bftheta_2)$, $\bftheta_1,\bftheta_2\!\in\!\RR^n$. Suppose $R\!\sim\!S^1(d,\RRR)$ and $\bfalpha\!=\!(\alpha_1,\dots,\alpha_n)\!\in\![0,\infty)^n$. Without loss of generality, assume $\alpha_1\!\le\!\dots\!\le\!\alpha_n$.
Denote the augmented process by $W:=( I\bfalpha, Y)$.
Proposition~\ref{propmultiparameter} states that $W(r)\!=\!(r\bfalpha,Y(r))\!\eqd\!(r\bfalpha,X(r\bfalpha))$ for $r\!\ge\!0$, thus proving the identity $I_1(\bftheta)\!=\!I_2(\bftheta)$, where~$I_1(\bftheta):=\int_{(0,\infty)}(\Phi_{(r\bfalpha,X(r\bfalpha))}(\bftheta)\!-\!1)\RRR(\rmd r)$ and $I_2(\bftheta):=\int_{(0,\infty)}(\Phi_{W(r(\bfe,\bfe))}(\bftheta)\!-\!1)\RRR(\rmd r)$.

Note $R\bfalpha\!\sim\!S^n(d\bfalpha,\RRR\circ (I\bfalpha)^{-1})$. Proposition~\ref{propcharexpo} and the transformation theorem tells us that $\Psi_{(R\bfalpha,X\odot(R\bfalpha))}(\bftheta)\!=\!\rmi d\skal{\bfalpha}{\bftheta_1}\!+\!d \Psi_Y(\bftheta_2)\!+\!I_1(\bftheta)$, also recalling $\bfalpha\tr\Psi_X(\bftheta_2)\!=\!\Psi_Y(\bftheta_2)$. Next, observe that~$R(\bfe,\bfe)\!\sim\!S^{2n}(d(\bfe,\bfe),\RRR\circ (I(\bfe,\bfe))^{-1})$ and $\Psi_{W\odot (R(\bfe,\bfe))}(\bftheta)\!=\!\Psi_{(R(\bfe,\bfe),W\odot (R(\bfe,\bfe)))}({\bf 0},\bftheta)$.
By Proposition~\ref{propcharexpo} and the transformation theorem, the RHS evaluates to
$d (\bfe,\bfe)\tr\Psi_W(\bftheta)\!+\!I_2(\bftheta)\!=\!\rmi d\skal{\bfalpha}{\bftheta_1}+d \Psi_Y(\bftheta_2)\!+\!I_2(\bftheta)$ by Proposition~\ref{propmultiparameter}.

The last statement in Proposition \ref{propweakalongray} follows from Proposition~\ref{propsupextendssub}.\halmos
\begin{remark}\label{remweakextstrsubord}
Let $B,B^*,N$ be independent processes, where $B\eqd B^*$ are standard Brownian motions, and $N$ is a Poisson process with unit rate. By Proposition~\ref{propweakalongray}, it follows from independence that $((I,2I),(B,B)\odot (I,2I)\big)
		\eqd(I,2I,B,B\!+\!B^*)\circ (I,I,I,I)$
and $((N,2N),(B,B)\odot (N,2N)\big)
		\eqd(I,2I,B,B\!+\!B^*)\circ (N,N,N,N)$.
Thus, we can represent these processes using strong subordination with the univariate subordinators $I$ and $N$, respectively.\halmos
\end{remark}
\noindent{\bf Moments.}~We give formulae for expected values and covariances.
\begin{proposition}\label{propmoments}
	If $X$ and $T$ be as in Definition~\ref{defsupervision}, then, for $t>0$,
\begin{align*}
\EE[T(t)]/t&=\bfd\!+\!\int_{[0,\infty)_*^n} \bft\TTT(\rmd\bft)\,,\quad\mbox{\em Cov}(T(t))/t=\int_{[0,\infty)_*^n} \bft'\bft\;\TTT(\rmd\bft)\,,\\
\EE[X\odot T(t)]/t&=\bfd\!\tr\!\bfmu+\int_{\DD^C}\! \bfx\,(\bfd\!\tr\!\XXX)(\rmd\bfx)+\int_{[0,\infty)_*^n}\!\EE[X(\bft)]\;\TTT(\rmd\bft)\,,\\
\mbox{\em Cov}(X\odot T(t))/t&=\bfd\tr\Sigma\!+\!\int_{\RR^n_*}\!\bfx'\bfx(\bfd\tr\XXX)(\rmd \bfx)\!+\!\int_{[0,\infty)_*^n}\!\!\EE[X'(\bft)X(\bft)]\TTT(\rmd\bft),\\
\mbox{\em Cov}(X\odot T(t),T(t))/t&=\int_{[0,\infty)_*^n} \EE[X'(\bft)]\bft\;\TTT(\rmd\bft)\,,
\end{align*}
provided the participating integrals are finite.
\end{proposition}
\noindent{\bf Proof.} Given the characteristics of $Z\eqd(T,X\odot T)\sim L^{2n}(\bfm,\Theta,\ZZZ)
$ in~\eqref{defmu1}--\eqref{PliY}, these follow from the general formulae for moments of L\'evy processes (see~\cite{s}, his Example~25.12).\halmos
\begin{remark}\label{remBMmoments} Let Brownian motion $B\sim BM^n(\bfmu,\Sigma)$ be the weak subordinate and $T\sim S^n(\bfd,\TTT)$ be the subordinator. By Proposition~\ref{propmoments}, for $1\!\le\!k\!\le\!n$,
\begin{eqnarray*}
\EE[(B\odot T)_k(1)]&=& \mu_k{\EE[T_k(1)]}\,,\\
\myVar((B\odot T)_k(1))&=&
\Sigma_{kk}\EE[T_k(1)]+\mu_k^2 \myVar(T_k(1))\,.
\end{eqnarray*}
Assume $1\!\le\!k\!\neq\!l\!\le\!n, u\!>\!0$, and set \[\tau_{k,l}(u)\,:=\,\TTT(\{\bft\!=\!(t_1,\dots,t_n)\!\in\![0,\infty)^n_*\!:\!t_k\!\wedge\!t_l\!>u\})\,.\] Recall~$s\wedge t\!=
\int_{(0,\infty)}\eins_{(u,\infty)}(s)\eins_{(u,\infty)}(t)\rmd u$, $s,t\!\ge\!0$, and $\int_{[0,\infty)_*^n}t_k\wedge t_l\TTT(\rmd\bft)\!=\!\int_{(0,\infty)}\tau_{k,l}(u)\rmd u$. Proposition~\ref{propmoments} states that
\[\mbox{Cov}((B\odot T)_k(1),(B\odot T)_l(1)))
\!=\!\mu_k\mu_l\mbox{Cov}(T_k(1),T_l(1))+\Sigma_{kl}(d_k\wedge d_l)+\Sigma_{kl}\int_{(0,\infty)}\tau_{k,l}(u)\;\rmd u\,.
\]
\halmos
\end{remark}
\begin{remark}\label{remmomentsVaG} Let $B\odot T\!\sim\!WV\bfalpha G^n(a,b,\bfalpha,\bfmu,\Sigma)$ be as in \eqref{defweakValphaG}. As the components of $T=(T_1,\dots,T_n)$ are standard gamma subordinators $T_k\sim\Gamma_S(b/\alpha_k)$, $1\!\le\! k\!\le\! n$, the first and second moments of an $\bfalpha G$-subordinator are determined as follows (see~\cite{Se08}),
\[\EE[T_k(1)]=1\,,\quad \myVar(T_k(1))=\alpha_k/b\/,\quad 1\le k\le n\,,\]
and
\[\myCov(T_k(1),T_l(1))=\alpha_k\alpha_l\myVar(T_0(1))=\alpha_k\alpha_la/b^2\,,\quad 1\le k\neq l\le n\,.\]
It follows from~Remark~\ref{remBMmoments} that, for $1\!\le\! k\!\le\! n$,
\begin{align*}
\EE[(B\odot T)_k(1)]=\mu_k, \;\;
\myVar((B\odot T)_k(1))=(b\Sigma_{kk}+\mu_k^2\alpha_k)/b\,,
\end{align*}
and these formulae match, not surprisingly, those of univariate $VG^1$-processes in~\cite{MaSe90} because of~\eqref{VGcomps}. If $1\!\le\! k\!\neq\! l\!\le\! n$, observe $\int_{(0,\infty)}\tau_{k,l;a,b,\bfalpha}(u)\;\rmd u=(\alpha_k\wedge \alpha_l) \EE[T_0(1)]=(\alpha_k\wedge \alpha_l) a/b$, with covariance given by
\begin{align}
\myCov((B\!\odot\!T)_k(1),(B\!\odot\!T)_l(1))&=\Sigma_{kl}
(\alpha_k\!\wedge\!\alpha_l)\EE[T_0(1)]\!+\!\mu_k\mu_l \myCov(T_k(1),T_l(1))\nonumber\\
&=(ab
(\alpha_k\wedge \alpha_l)\,\Sigma_{kl}
+a\alpha_k\alpha_l\mu_k\mu_l)/b^2\,.\label{eqaddcorr}
\end{align}

These moments for the $WV\bfalpha G^n$-process have also been derived in \cite{MiSz17} as well as higher moments.

For traditional subordination, \eqref{eqaddcorr} reduces to $\alpha_k\alpha_l\mu_k\mu_l/b^2$~(see~\cite{Se08}, her Section~4) as $\Sigma$ is diagonal, which was noted as a disadvantage in~\cite{Gu13,LS10}.
In contrast, $B\odot T$ has an additional cor\-relation term which includes the cor\-relation of the Brownian motion.\halmos\end{remark}
\noindent{\bf Superposition.}~If a process $X$ is weakly subordinated by a superposition of several independent subordinators, then its law equals the sum of independent L\'evy processes.
\begin{proposition}\label{propsupsubord} Let $X$ be an $n$-dimensional L\'evy process.
Let $\bfd\in[0,\infty)^n$ be a deterministic vector. If~$T_1,\dots,\allowbreak T_m$ are independent $n$-dimensional drift-less subordinators, then~$T\!:=\!I \bfd+\sum_{k=1}^mT_k$
is an $n$-dimensional subordinator with drift $\bfd$ and $(T,X\odot T)\eqd\sum_{j=0}^mA_j$, where $A_0,A_1,\dots,\allowbreak A_m$ are independent L\'evy processes with $A_0\eqd (I\bfd,X\odot I\bfd)$, $A_k\eqd (T_k,X\odot T_k)$, $1\!\le\! k\!\le\! m$.
\end{proposition}
\noindent{\bf Proof.}\label{proofpropsupsubord}~Assume that $T_1,\dots,T_m,A_0,\dots,A_m$ are independent processes, where~$T_k\!\sim\! S^n({\bf 0},\TTT_k)$, $A_0\!\eqd\!(I\bfd,X\odot I\bfd)$ and
$A_k\!\eqd\! (T_k,X\odot T_k)$, $1\!\le\!k\!\le\!m$. In particular, note $T\sim S^n({\bf d},\sum_{k=1}^m \TTT_k)$, then by~\eqref{propsupervisexpo},
\begin{eqnarray*}\Psi_{(T,X\odot T)}(\bftheta)&=&
\rmi\skal{\bftheta_1}{\bfd}+(\bfd\tr\Psi_{X})(\bftheta_2)+\int_{[0,\infty)^n_*}
(\Phi_{(\bft,X(\bft))}(\bftheta)-1)\;(\sum_{k=1}^m \TTT_k)(\rmd \bft)\\
&=&\sum_{k=0}^m \Psi_{A_k}(\bftheta)=\Psi_{\sum_{k=0}^m A_k}(\bftheta)\,,\quad
\bftheta=(\bftheta_1,\bftheta_2),\quad \bftheta_1,\bftheta_2\in\RR^n\,,
\end{eqnarray*}
as desired.\halmos
\begin{remark}\label{remexsupsubord}
	In the context of traditional subordination (see~\cite{BKMS16}, their Proposition~4.1), Proposition~\ref{propsupsubord} holds without assuming drift-less subordinators.
	This is more delicate when dealing with weak subordination.
	Let $B,B^*,W,W^*$ be independent standard univariate Brownian motions.

Remark~\ref{rem2counterexamplBBT2T} states that $(B,B)\odot (I,2I)\eqd(B,B\!+\!B^*)$ and $(B,B)\odot (2I,I)\eqd (W\!+\!W^*,W)$. Proposition~\ref{propsupextendssub} states that $(B,B)\odot(3I,3I)\eqd(B,B)\circ(3I,3I)$. Note $(B,B+B^*)+(W+W^*,\allowbreak W)\sim BM^2({\bf 0},[(3,2),(2,3)])$ and
	$(B,B)\circ(3I,3I)\!\sim\! BM^2({\bf 0},[(3,3),(3,3)])$.
	There are no independent processes $Y_1,Y_2$ such that, simultaneously,
	$Y_1\eqd (B,B)\odot(I,2I)$, $Y_2\eqd (B,B)\odot (2I,I))$ and $Y_1\!+\!Y_2\eqd(B,B)\odot(3I,3I)$.
	\halmos
\end{remark}
\begin{remark}\label{remVaGdecomposition} Let $B\odot T\sim WV\bfalpha G^n(a,b,\bfalpha,\bfmu,\Sigma)$ in~\eqref{defweakValphaG}. We derive a joint representation of $(T,B\odot T)$ in terms of a superposition of gamma processes and variance-gamma processes. Let $B,B^{(1)},\dots,B^{(n)},W^{(\bfalpha)}, G_0,\dots,\allowbreak G_n$ be independent, where $B^{(1)},\dots,B^{(n)}$ are copies of $B\sim BM^n(\bfmu,\Sigma)$, $G_0,\dots,G_n$ are as in~\eqref{semeraro2} and $W^{(\bfalpha)}\sim BM^n(\bfalpha\tr\bfmu,\bfalpha\tr\Sigma)$
is a Brownian motion.

Next, standardise $bG_0/a\sim\Gamma_S(a)$ and $(b/(b\!-\!a\alpha_k))G_k\sim \Gamma_S(\beta_k)$
to see that $V_0\!:=\!W^{(\bfalpha)}\!\circ\!(G_0\bfe)\!\sim\!VG^n(a,(a/b)(\bfalpha\tr\bfmu,\bfalpha\tr\Sigma))$
and $V_k:=B^{(k)}_k\!\circ\!G_k\!\sim\!VG^1(\beta_k,((b\!-\!a\alpha_k)/b)\,(\mu_k,\Sigma_{kk}))$, $1\!\le\! k\!\le\! n$. Note $V_0,\dots,V_n$ are independent.

Plainly, $T$ in~\eqref{semeraro2} is the superposition
of independent univariate gamma processes travelling along deterministic rays generated by  $\bfalpha,\bfe_1,\dots,\bfe_n\in[0,\infty)^n_*$. Combining Propositions~\ref{propsupsubord} and~\ref{propweakalongray} yields, for $Z\eqd (T,B\odot T)$,
\begin{eqnarray*}
\lefteqn{Z\,\eqd\,
	(G_0\bfalpha,B\odot(G_0\bfalpha))+\sum_{k=1}^n(G_k\bfe_k,B^{(k)}\odot(G_k\bfe_k))}&&\\
&\eqd&\!\!(I\bfalpha,W^{(\bfalpha)})\!\circ\!(G_0(\bfe,\bfe))+\sum_{k=1}^n(G_k\bfe_k,(B^{(k)}_k\circ G_k)\bfe_k)\!=\!(G_0\bfalpha,V_0)+\sum_{k=1}^n(G_k\bfe_k,V_k\bfe_k)\,.
\end{eqnarray*}
Our~$WV\bfalpha G$-process satisfies $Y\eqd B\odot T\eqd V_0+\sum_{k=1}^nV_k\bfe_k$ as the superposition of independent $VG^n$-processes~(for the strong formulation, see~\cite{BKMS16}, their Remark~2.17).\halmos
\end{remark}
\noindent{\bf Subordinators with independent components.}~If a drift-less subordinator has independent components, then so does any associated weakly subordinated process.
\begin{proposition}\label{propsubordindependent}Let $X$ and $T$ be as in Definition~\ref{defsupervision}, with
drift-less $T$. If the components of $T$ are independent, then so are those of $X\odot T$.
\end{proposition}
\noindent{\bf Proof.}~If $T=(T_1,\dots,T_n)\sim S^n({\bf 0},\TTT)$ has
independent components $T_1\sim S^1(0,\TTT_1),\dots,\allowbreak T_n\sim S^1(0,\TTT_n)$, then $\TTT=\sum_{k=1}^n\bfdelta_0^{\otimes(k\!-\!1)}\otimes \TTT_k\otimes \bfdelta_0^{\otimes(n\!-\!k)}$.
In~\eqref{SigmaY2}--\eqref{PliY2}, note $\bfd\tr\Sigma=0$ and $\ZZZ_2=\sum_{k=1}^n\bfdelta_0^{\otimes(k\!-\!1)}\otimes \YYY_k\otimes \bfdelta_0^{\otimes(n\!-\!k)}$,
where $\ZZZ_2,\YYY_1,\dots,\YYY_n$ are
the L\'evy measures corresponding to $X\odot T,X_1\odot T_1,\dots,X_n\odot T_n$, as required.\halmos
\begin{remark}
If $B,B^*,N,N^*$ are independent processes, such that $B,B^*$ are univariate standard Brownian motions
and $N,N^*$ are Poisson processes with unit rate, then it is straightforwardly
verified from \eqref{muY2}--\eqref{PliY2} and Proposition~\ref{propsupextendssub} that $(B,B)\odot (N,N^*)\eqd(B\circ N,B^*\circ N^*)$ decomposes into a L\'evy process with independent compound Poisson components.
\halmos
\end{remark}
We have previously listed sufficient conditions for strong subordination~\cite{BPS01,s} to stay in the class of L\'evy processes. Next, we show that these conditions are necessary in some cases (see Subsection~\ref{subsecproofIII} for a proof).
\begin{proposition}\label{propindepindist}
Let~$T=(T_1,T_2)$ and $X=(X_1,X_2)$ be independent bivariate L\'evy processes, where $T$ is a subordinator. Suppose neither $T_1\equiv 0$ nor $T_2\equiv 0$. If $X\circ T$ is also a L\'evy process, then $T_1=T_2$ must be indistinguishable, provided one of the following holds in addition:\\
	(i)~$X\eqd -X$ is symmetric, and $X_1,X_2$ are dependent;\\
	(ii)~$T$ is deterministic, and $X_1,X_2$ are dependent;\\
	(iii)~$T$ admits a finite first moment, and $X$ admits a finite second moment with correlated components $X_1,X_2$.\label{propcounter1}
\end{proposition}
In Proposition~\ref{propmonocase}, we stated monotonicity as a sufficient condition ensuring that
the weakly subordinated process matches the marginal distributions of the strongly subordinated one. Next,
we show that for this purpose, monotonicity is needed in some cases (see Subsection~\ref{subsecproofIII} for a proof).
\begin{proposition}\label{propmonotindepmarginal}
If~$T=(T_1,T_2),X=(X_1,X_2)$ and $Y=(Y_1,Y_2)$ are bivariate L\'evy processes, where $T,X$ are independent and $T$ is a subordinator,
while $X$ has dependent component $X_1,X_2$, then there is at least one $t\in(0,\infty)$ violating $X(T(t))\eqd Y(t)$, provided one of following holds in addition:\\
(i)~both $T,X$ admit finite second moments, while $X$ has correlated components and $T$ has non-monotonic components;\\
(ii)~$Y\eqd X\odot T$, $X$ is symmetric, while $T_1,T_2$ are independent, drift-less and nontrivial subordinators.
\end{proposition}
\begin{remark} In Proposition~\ref{propmonotindepmarginal}(ii), the subordinator has independent and non-deterministic components, and so is non-monotonic, that is, neither $T_1-T_2$ nor $T_1-T_2$ is a
subordinator.

It would be interesting to see whether or not the conditions in Propositions~\ref{propindepindist}--~\ref{propmonotindepmarginal}
could be further weakened. We speculate that this extension is possible based on Dynkin-type formulae and fluctuation theory for L\'evy processes. We have to leave this as an interesting avenue of future research.\halmos
\end{remark}
\section{Variance Generalised Gamma Convolutions}\label{secVGGC}
In this section, the weak subordinate is Brownian motion $B\sim BM^n(\bfmu,\Sigma)$, and
$T\sim S^n(\bfd,\TTT)$ is the subordinator. Since the L\'evy measure of $B$ is 0, we get simplifications in~\eqref{VVVsubord}--\eqref{PliY} and~\eqref{muY2}--\eqref{PliY2}.
The weakly subordinated process is denoted by $Y\eqd B\odot T\sim L^n(\bfm_2,\bfd\tr\Sigma,\YYY)$.\\[1mm]
Thorin~\cite{Th77a,Th77b} characterised the class of generalised gamma convolutions ($GGC$) as the subset of univariate Borel probability measures containing arbitrary finite convolutions of gamma distributions, while being closed under convergence in distribution (see the survey article~\cite{JLY08} and the monograph~\cite{SSV}).
Multivariate extensions of these results and examples
have been investigated in~\cite{BMS06,Bo09,PS14}, and these are subclasses
of the self-decomposable and, thus, infinitely divisible distributions. Our subordinators will be taken from this class.

\noindent{\bf Thorin subordinator.}~In our exposition we follow~\cite{BKMS16}. Recall $\ln^-\!x\!=\!-\!\eins_{(0,1]}\!(x)\!\ln\!x,x\!>\!0$. A nonnegative Borel measure $\UUU$ on $[0,\infty)^n_*$ is called an $n$-dimensional {\em Thorin measure}, provided
\begin{equation*}
\int_{[0,\infty)_*^n}\;\left(1\!+\!\ln^-\|\bfu\|\right)\wedge\left(1\big/\|\bfu\|\right)\;\UUU(\rmd \bfu)\quad <\quad \infty\,.
\end{equation*}
If $\bfd\in [0,\infty)^n$ and $\UUU$ is a Thorin measure, we call an $n$-dimensional subordinator $T$ a {\em Thorin subordinator}, in brief $T\sim\GGC_S^n(\bfd,\UUU)$, whenever, for all $t\ge0,\bflambda\in[0,\infty)^n$, it has Laplace exponent
\begin{equation}\label{LaplaceGGCd}
-\ln \EE\exp\{-\skal{\bflambda}{T(t)}\}\,=\,t\skal \bfd\bflambda+t\int_{[0,\infty)^n_*} \!\ln\big\{(\|\bfu\|^2+\skal{\bflambda}{\bfu})\big/\|\bfu\|^2\big\}\,\UUU(\rmd \bfu)\,.
\end{equation}
The distribution of a Thorin subordinator is uniquely determined by~$\bfd$ and~$\UUU$.

Let $\mySS_+\!:=\!\mySS\!\cap\![0,\infty)^n_*$, where $\mySS\!:=\!\{\bfs\!\in\!\RR^n\!:\!\|\bfs\|\!=\!1\}$ is the unit sphere. If $T\sim S^n(\bfd,\TTT)$, the L\'evy measure $\TTT$ is derived using a polar-decomposition
of its Thorin measure. Specifically, if $A\in[0,\infty)^n_*$ is a Borel set, then we may write (see~\cite{BKMS16}, their Lemma~4.1)
\[\UUU(A)\,=\,(\SSS\otimes\KKK)\circ \big((\bfs,r)\mapsto r\bfs\big)^{-1}(A)\,=\,\int_{\mySS_+}\int_{(0,\infty)} \KKK(\bfs,\rmd r)\eins_A(r\bfs)\,\SSS(\rmd\bfs)\,.\]
Here, $\SSS$ is a finite nonnegative Borel measure on $\mySS_+$ and $\KKK$ is a Thorin kernel, that is a nonnegative Borel kernel with
\begin{equation*}
0\,<\,\int_{(0,\infty)}
(1\!+\!\ln^- r)\wedge (1/r)\,\KKK(\bfs, \rmd r)\,<\,\infty\,,\qquad \bfs\in\mySS_+\,.\end{equation*}
Recall $\GGG_b$ is the L\'evy measure of a standard gamma subordinator with shape parameter~$b$.
\begin{lemma}
If $T\sim\GGC_S^n(\bfd,\UUU)$, then $T\sim S^n(\bfd,\TTT)$, where
\begin{equation}\label{GGGClevmeasintermsstandardGamma}
\TTT\,=\,\Big\{
\frac{\UUU(\rmd \bfu)}{\|\bfu\|^2}\otimes\GGG_{\|\bfu\|^2}(\rmd g)\Big\}\circ \big((\bfu,g)\mapsto g\bfu\big)^{-1}\,.\end{equation}
\end{lemma}
\noindent{\bf Proof.}~If $T\sim\GGC_S^n(\bfd,\UUU)$, then
$T\sim S^n(\bfd,\TTT)$ in polar coordinates is (see~\cite{BKMS16}, their Equations (2.17)--(2.18))
\begin{align}\label{GGClevydensity}
\TTT(A)&=\int_{\mySS^+}\int_{(0,\infty)}\eins_A(r\bfs)k(\bfs,r)\,\frac{\rmd r}r\,\SSS(\rmd \bfs)\,,\qquad \mbox{$A\subseteq [0,\infty)^n_*$ Borel}\,,\\
k(\bfs,r)&=\int_{(0,\infty)} e^{-r v}\,\KKK(\bfs,\rmd v)\,,\;\quad r>0,\bfs\in\mySS_+\,.\nonumber
\end{align}
If $A\subseteq [0,\infty)^n$ is Borel, by using~\eqref{GGClevydensity} and making the substitution $g=r/\|\bfu\|$, we get that
\[\TTT(A)\,\,=\int_{[0,\infty)_*^n}\left\{\|\bfu\|^2\int_{(0,\infty)}\eins_{A}(g\bfu)\;e^{-g\|\bfu\|^2}\,\frac{\rmd g}g\right\}\frac{\UUU(\rmd\bfu)}{\|\bfu\|^2}\,.\]
Here, the RHS matches the RHS of~\eqref{GGGClevmeasintermsstandardGamma} when evaluated at $A$.\halmos
\begin{remark} \label{remTabalphaGGC}
If $T\sim \bfalpha G^n(a,b,\bfalpha)$ in \eqref{semeraro2}, then $T$ is determined as the superposition of independent gamma subordinators $G_0,\dots,G_n$, travelling along rays generated by $\bfalpha,\bfe_1,\dots,\bfe_n$, respectively.
Recall $\beta_k:=(b\!-\!a\alpha_k)/\alpha_k$, $1\!\le\!k\!\le\!n$. If $\bflambda\in[0,\infty)^n_*$, then
\begin{eqnarray*}
\lefteqn{-\ln\EE\exp\{-\skal\bflambda{T(t)}\}}&&\\
&=&-\ln \EE[\exp\{-G_0(t)\skal\bflambda\bfalpha\}]-\sum_{k=1}^n\ln\EE[\exp\{-G_k(t)\skal\bflambda{\bfe_k}\}]\\
&=&at\ln\{(b+\skal \bflambda\bfalpha)/b\}+\sum_{k=1}^n\beta_k t\ln\{((b/\alpha_k)+\skal\bflambda{\bfe_k})/(b/\alpha_k)\}\,.
\end{eqnarray*}
Here, we used independence and the Laplace exponent of the underlying gamma subordinators.

The RHS matches~\eqref{LaplaceGGCd} for $\bfd\!=\!{\bf 0}$ and $\UUU_{a,b,\bfalpha}:=a\bfdelta_{b\bfalpha/\|\bfalpha\|^2}\!+\!\sum_{k=1}^n\beta_k\,\bfdelta_{b\bfe_k/\alpha_k}$.
Therefore,~$\UUU_{a,b,\bfalpha}$ defines a finitely supported Thorin measure, and $T\sim\GGC^n_S({\bf 0},\UUU_{a,b,\bfalpha})$ is a drift-less Thorin subordinator.

Using $\UUU_{a,b,\bfalpha}$ and~\eqref{GGGClevmeasintermsstandardGamma}, it is possible
to give an alternative derivation of the
L\'evy measure $\TTT_{a,b,\bfalpha}$ in~\eqref{levyalphasubord} (see~\cite{BKMS16} and~\cite{LS10}, their Lemma~2.13 and  their Theorem~1.1, respectively).\halmos
\end{remark}
\noindent{\bf Variance generalised gamma convolutions.}~For the parameters of this model we assume an $n$-dimensional Thorin measure $\UUU$, $\bfmu\in\RR^n$, $\bfd\in[0,\infty)^n$ and a covariance matrix $\Sigma\in\RR^{n\times n}$.
Let $B\sim BM^n(\bfmu,\Sigma)$ be a Brownian motion.
Let $T\sim \GGC_S^n(\bfd,\UUU)$. Given such $B$ and $T$, we call a L\'evy process of the form $Y\eqd B\odot T$ an $n$-dimensional {\em variance generalised gamma convolution $($$\VGGC^{n}$$)$} process with parameters $\bfd,\bfmu,\Sigma,\UUU$. We write this as
\begin{equation*}
Y\sim \VGGC^{n}(\bfd,\bfmu,\Sigma,\UUU):=BM^n(\bfmu,\Sigma)\odot \GGC^n_S(\bfd,\UUU)\,.\end{equation*}
Theorem~\ref{thmexist} ensures the existence of $Y\sim\VGGC^{n}(\bfd,\bfmu,\Sigma,\UUU)$.\\[1mm]
\noindent{\bf Characteristics.}~We derive formulae of the characteristic exponent and the L\'evy measure, valid within the $\VGGC^n$-class.
If $\emptyset\!\neq\!J\!\subseteq\!\{1,\dots,n\}$, introduce $C_J\!\subseteq\!V_J\!\subseteq\!\RR^n$, where
$\bfu\!=\!(u_1,\dots,u_n)\!\in\!C_J$ and $\bfy=(y_1,\dots,y_n)\!\in\!V_J$ if and only if $u_j\!>\!0$ for all $j\!\in\!J$ and $y_j\!\neq\!0$
for all $j\!\in\!J$, respectively. If $\bfu\in C_J$, while $\Sigma$ is invertible, the restriction $(\bfu\tr\Sigma)_J:\RR^n\bfpi_J\to\RR^n\bfpi_J$,
$\bfx\mapsto\bfx(\bfu\tr\Sigma)_J:=\bfx(\bfu\tr\Sigma)$ is an invertible linear mapping, thus having inverse $(\bfu\tr\Sigma)^{-1}_J$ and determinant $|\bfu\tr\Sigma|_J$.
\begin{theorem}\label{thmGGCcharfct} If $Y\sim \VGGC^{n}(\bfd,\bfmu,\Sigma,\UUU)$, then $Y\sim L^{n}(\bfm_2,\bfd\tr\Sigma,\YYY)$, where $\bfm_2=\bfd\tr\bfmu+\int_{\DD_*}\bfy\,\YYY(\rmd\bfy)$, $\VVV$ is the quantity in~\eqref{levyVGgen} and
\begin{equation}\label{levyVVGGCintermsofVGg}
\YYY\,=\,\Big\{
\frac{\UUU(\rmd \bfu)}{\|\bfu\|^2}\otimes\VVV_{\|\bfu\|^2,\bfu\tr\bfmu,\bfu\tr\Sigma}(\rmd\bfy)
\Big\}\circ \big((\bfu,\bfy)\mapsto \bfy\big)^{-1}\,,\end{equation}
and,~for $\bftheta\!\in\!\RR^n$,
\begin{equation}	\label{GVGcharexpo}
\Psi_Y(\bftheta)\,=\,
\rmi \skal{\bfd\tr\bfmu}{\bftheta}-\frac 12 \|\bftheta\|^2_{\bfd\tr\Sigma}
	-\int_{[0,\infty)^n_*}\!\ln\big\{(\|\bfu\|^2
-\rmi\skal{\bfu\tr\bfmu}{\bftheta}+\frac 12\|\bftheta\|^2_{\bfu\tr\Sigma})\big/\|\bfu\|^2\big\}\UUU(\rmd\bfu).
\end{equation}
If, in addition, $\Sigma$ is invertible, then $\YYY=\sum_{\emptyset\neq J\subseteq\{1,\dots,n\}}\YYY_J$, $\YYY_J(\RR_J\backslash V_J)=0$, where $\YYY_J$ is absolutely continuous with respect to $\rmd\bfy\circ\bfpi^{-1}_J$ having density $v_J(\bfy)
=\int_{C_J}\,\nu_{J}(\bfu,\bfy)\UUU(\rmd\bfu)$, where $\bfu\!\in\!C_J$, $\bfy\!\in\!V_J$, $c_J:=2/(2\pi)^{\#J/2}$, and \begin{eqnarray*}
\nu_J(\bfu,\bfy)&=&c_J\KKKK_{\#J/2}\{[\|\bfy\|_{(\bfu\tr\Sigma)^{-1}_J}(2\|\bfu\|^2\!+\!\|\bfu\tr\bfmu\|^2_{(\bfu\tr\Sigma)_J^{-1}})^{1/2}\}\label{defhJ}\\
&&{}\times\exp\big\{\skal{\bfy}{\bfu\tr\bfmu}_{(\bfu\tr\Sigma)_J^{-1}}\big\}\Big/
\big\{|\bfu\tr\Sigma|_J^{1/2}\,\|\bfy\|^{n}_{(\bfu\tr\Sigma)_J^{-1}}\big\}\,.\nonumber
\end{eqnarray*}
\end{theorem}
\noindent{\bf Proof.}~The formulae of the triplet $(\bfm_2,\bfd\tr\Sigma,\YYY)$ follow from Proposition~\ref{propprops}. To see this, let $A\subseteq\RR^n_*$ be a Borel set. Combining~\eqref{levyVGgen} with Proposition~\ref{propmultiparameter} yields $\int_{(0,\infty)}\PP(B(g\bfu)\in A)\;\GGG_{\|\bfu\|^2}(\rmd g)=\VVV_{\|\bfu\|^2,\bfu\tr\bfmu,\bfu\tr\Sigma}(A)$. In particular, we get from~\eqref{PliY2} and~\eqref{GGGClevmeasintermsstandardGamma} that
\[\YYY(A)\,=\,\int_{[0,\infty)^n_*}\int_{(0,\infty)}\PP(B(g\bfu)\in A)\;\GGG_{\|\bfu\|^2}
(\rmd g)\,\frac{\UUU(\rmd\bfu)}{\|\bfu\|^2}\,,\]
where the RHS matches the RHS in~\eqref{levyVVGGCintermsofVGg} when evaluated at $A$.

As $\bft\mapsto\EE[B(\bft)\eins_{\DD_*}(B(\bft))]$ is $\TTT$-integrable by \eqref{PiFVintegrab}~and~\eqref{m1truncabsoutineq}, $\bfy\mapsto \bfy\eins_{\DD_*}(\bfy)$ is
$\YYY(\rmd\bfy)=\PP(B(\bft)\in \rmd\bfy)\TTT(\rmd\bft)$-integrable by the transformation theorem.
In particular, the linear term under the integral in \eqref{0.1} cancels, and
combining~\eqref{charexpoVG} and~\eqref{levyVVGGCintermsofVGg} to see that $\int_{\RR^n_*}e^{\rmi\skal\bftheta\bfy}-1\,\YYY(\rmd\bfy)$ matches the integral in
\eqref{GVGcharexpo}.

In view of~\eqref{defgVGbmuSigma} and~\eqref{levyVVGGCintermsofVGg}, the L\'evy density formula follows straightforwardly.\halmos
\begin{remark}\label{remVGGn1VGGnnsubclVGGn}~Strong univariate subordination
of an arbitrary Brownian motion with an independent univariate Thorin subordinator was investigated in~\cite{Gr07}.
The corresponding class of L\'evy processes was called $\VGGC^{n,1}$ in~\cite{BKMS16}.
Using our notation, we have  $\VGGC^{n,1}(d,\bfmu,\Sigma,\UUU_0):=
\VGGC^n(d \bfe,\bfmu,\Sigma,\int_{(0,\infty)}\bfdelta_{u\bfe}\,\UUU_0(\rmd u))$,
where $\bfmu\in\RR^n$, $d\in [0,\infty)$,
while $\Sigma\in\RR^{n\times n}$ is an arbitrary covariance matrix and $\UUU_0$ is a univariate Thorin measure.
The $VG^n$-process~\cite{MaSe90} provides us with an example of a
$\VGGC^{n,1}$-process.

The $\VGGC^{n,n}$-class was introduced in~\cite{BKMS16} to complement the $\VGGC^{n,1}$-class and contains processes formed by strong multivariate subordination of an {\em independent-component} Brownian motion with a Thorin subordinator. More specifically, $\VGGC^{n,n}(\bfd,\bfmu,\Sigma,\UUU):=\VGGC^{n}(\bfd,\bfmu,\Sigma,\UUU)$
where $\bfd\in[0,\infty)^n$, $\bfmu\in\RR^n$, while $\Sigma$ is a
covariance matrix of {\em diagonal} form and $\UUU$ is an $n$-dimensional Thorin measure. The strong $V\bfalpha G$-process~\cite{Se08} is an example of a $\VGGC^{n,n}$-process.

In~\cite{BKMS16} (see Part~(i) of their Theorems~2.3 and 2.5), formulae of the characteristic exponents of $\VGGC^{n,1}\cup \VGGC^{n,n}$-processes are stated separately, while our Theorem~\ref{thmGGCcharfct} unifies both classes as special cases.\halmos
\end{remark}
\begin{remark}\label{rem1VaGVGGC}
Though it does not need to be an element of the $\VGGC^{n,1}\cup\VGGC^{n,n}$-class, a $WV\bfalpha G$-process always belongs to the $\VGGC^n$-class.

If $Y\eqd B\odot T\sim WV\bfalpha G^n(a,b,\bfalpha,\bfmu,\Sigma)$ in~\eqref{defweakValphaG}, then $Y$ is also a $\VGGC^n$-process. More specifically, with $\UUU_{a,b,\bfalpha}$ as in
Remark~\ref{remTabalphaGGC},
we have
 \[Y\eqd B\odot T\sim BM^n(\bfmu,\Sigma)\odot \GGC^n_S({\bf0},\UUU_{a,b,\bfalpha})=\VGGC^{n}({\bf 0},\bfmu,\Sigma,\UUU_{a,b,\bfalpha})\,.\]
In particular,
it follows from Theorem~\ref{thmGGCcharfct} that, $\bftheta=(\theta_1,\dots,\theta_n)\in\RR^n$,
\begin{eqnarray} \nonumber
&&\Psi_Y(\bftheta)=\Psi_{B\odot T}(\bftheta)\\&=&
-a\;\ln\big\{(b\!-\!\rmi \skal{\bfalpha\tr\bfmu}\bftheta+\frac 12\|\bftheta\|^2_{\bfalpha\tr\Sigma})\big/b\big\}\label{v-alpha-g-cf}-\sum_{k=1}^n\beta_k\ln\big\{(b\!-\!\rmi \alpha_k\mu_k\theta_k+\frac 12\alpha_k\theta_k^2\Sigma_{kk})\big/b\big\}\,.\nonumber
\end{eqnarray}
If, in addition, $\Sigma$ is invertible, then combining Theorem~\ref{thmGGCcharfct} and \eqref{defgVGbmuSigma}, and recalling $K_{1/2}(r)=\pi^{1/2}e^{-r}(2r)^{-1/2}$, $r>0$, (see~\cite{GrRy96}, their Equation (8.469)--3), we find that $Y$ has L\'evy measure $\YYY$ given by
\[\YYY(\rmd \bfy)={f_0(\bf y)\rmd \bfy}+ \sum_{k=1}^n \bfdelta^{\otimes(k-1)}_0 \otimes (f_k(y_k) \rmd y_k)\otimes \bfdelta^{\otimes(n-k)}_0\,,\]
where $\bfy=(y_1,\dots,y_n)\in\RR^n_*$,
\[f_0(\bfy)\,=\,\frac{2a\exp\big\{\skal{\bfy}{\bfalpha\!\tr\!\bfmu}_{(\bfalpha\tr\Sigma)^{-1}}\big\}}{(2\pi)^{n/2}|\bfalpha\!\tr\!\Sigma|^{1/2}
	\|\bfy\|^n_{(\bfalpha\tr\Sigma)^{-1}}}\;\;\KKKK_{n/2}\big\{\|\bfy\|_{(\bfalpha\tr\Sigma)^{-1}}(2b+\|\bfalpha\tr\bfmu\|^2_{(\bfalpha\tr\Sigma)^{-1}})^{1/2}\}\,,\]
and $1\le k\le n$, $y\in\RR_*$,
\[f_k(y)\,=\,\frac{\beta_k}{|y|}\exp\left\{\big(
\alpha_k^{1/2}\mu_k y-|y|(2b \Sigma_{kk}+\alpha_k\mu_k^2)^{1/2}\big)\Big/(\alpha_k^{1/2}\Sigma_{kk})
\right\}.\]
Alternatively, this decomposition could be derived from Remark~\ref{remVaGdecomposition}.\halmos
\end{remark}
\noindent{\bf Sample Paths.}~To see how sample path properties such as $q$-variation of the Thorin subordinator is propagated through Brownian motion, we generalise the corresponding result in~\cite{BKMS16} (see their Propositions~2.1--2.2; see Subsection~\ref{subsecproofIV} for a proof).
\begin{proposition}\label{propVGGCFV}\hspace*{-1mm}Let $T\sim \GGC_S^n(\bfd,\UUU)$ and $Y\sim VGG^n(\bfd,\bfmu,\Sigma,\UUU)$ with L\'evy measures $\TTT$ and $\YYY$, respectively. Suppose $0\!<\!q\!<\!1$.\\
(i)~$\int_{\DD_*^C}\UUU(\rmd\bfu)/\|\bfu\|^{q}$ is finite if and only if $\int_{\DD_*} \|\bft\|^q\,\TTT(\rmd\bft)$ is.\\
(ii)~If  $\int_{\DD^C}\UUU(\rmd\bfu)/\|\bfu\|^{q}$ is finite, then $\int_{\DD_*}\|\bfy\|^{2q}\,\YYY(\rmd\bfy)$ is finite. If $\Sigma$ is invertible, then also the converse holds.\\
(iii)~If $\bfd={\bf 0}$ and $\int_{\DD^C}\UUU(\rmd\bfu)/\|\bfu\|^{1/2}$ is finite, then $Y$ is a drift-less $FV^n$-process. If, in addition, $\Sigma$ is invertible, then
$Y\!\sim\!FV^n$ implies~$\bfd\!=\!{\bf 0}$ and $\int_{\DD^C}\UUU(\rmd\bfu)/\|\bfu\|^{1/2}<\infty$.
\end{proposition}
\begin{remark} In~\cite{BKMS16} (see their Remarks~2.8--2.9), examples are found of drift-less multivariate Thorin subordinators subordinating Brownian motion in the strong sense,
with the associated
$\VGGC^{n,1}\cup\VGGC^{n,n}$-process having sample paths of unbounded variation.
Proposition~\ref{propVGGCFV} states that those examples have counterparts in the weak sense.\halmos
\end{remark}
\begin{remark} If $\UUU$ is a finitely supported nonnegative measure on $[0,\infty)^n_*$, then $\UUU$ is in particular
a Thorin measure, and any associated drift-less $\VGGC^n$-process must be a $FV^n$-process as is straightforwardly derived from Proposition~\ref{propVGGCFV}(iii).
In particular, see~Remark~\ref{remTabalphaGGC}, weak variance-$\bfalpha$-gamma processes are drift-less $FV^n$-processes.\halmos
\end{remark}
\begin{remark} Weak subordination has applications in financial modelling. In \cite{MiSz17}, log returns of multiple dependent prices and $VG$-marginals were modelled using a $WV\bfalpha G^n$-process. In \cite{Mi17}, the log returns were modelled using a weakly subordinated process where the subordinator, interpreted as an information flow process, has jump dependence specified by a L\'evy copula while its marginals may be chosen arbitrarily. In \cite{Ma17}, $WV\bfalpha G^n$-processes were applied to instantaneous portfolio theory. In our future work, we will consider statistical inference for $WV\bfalpha G^n$-processes~\cite{BLM17b}, and conditions for the self-decomposability within the weak $VGG^n$-class~\cite{BLM17a}. \halmos
\end{remark}
\section{Proofs}\label{secproofs}
\subsection{Proof of Proposition~\ref{propmultiparameter}}\label{subsecproofI}
For $\bft=(t_1,\dots,t_n)\in[0,\infty)^n$ recall that $\langle(1),\dots,(n)\rangle$ denotes the associated permutation of the
ordering $t_{(1)}\le{\dots}\le t_{(n)}$ and $\Delta t_{(k)}$ correspond to its $k$th spacing. For $1\!\le\!m\!\le\!n$, let $\bfpi_m:=\bfpi_{\{(m),\dots,(n)\}}:\RR^n\to\RR^n$, $\bftheta\mapsto\bftheta\bfpi_m$.
Let $X\sim L^n(\bfmu,\Sigma,\XXX)$ with $\Psi$ as in~\eqref{0.1}. For $\bftheta=(\theta_1,\dots,\theta_n)\in\RR^n$, we have
\[
\sum_{k=1}^n\theta_k X_k(t_k)=\sum_{k=1}^n\theta_{(k)} X_{(k)}(t_{(k)})=\sum_{k=1}^n\sum_{m=1}^k \theta_{(k)}\big(X_{(k)}(t_{(m)})\!-\!X_{(k)}(t_{(m-1)})\big)\,,
\]
and thus, by interchanging the order of summation on the RHS,
\[
\sum_{k=1}^n\theta_k X_k(t_k)=\sum_{m=1}^n\sum_{k=m}^n\theta_{(k)}\Big(X_{(k)}(t_{(m)})\!-\!X_{(k)}(t_{(m-1)})\Big)\,,
\]
giving $\EE\exp(\rmi\skal\bftheta{X(\bft)})=\exp\{\sum_{m=1}^n \Delta t_{(m)}\Psi(\bftheta\bfpi_m)\}$
which matches~\eqref{multiexponent1}.

Since projections are self-adjoint, we must have
\[
\Psi(\bftheta\bfpi_m)=
\rmi \skal {\bfmu\bfpi_m}\bftheta-\frac 12\;\|\bftheta\bfpi_m\|^2_{\Sigma}+\int_{\RR_*^n}\big(e^{\rmi\skal\bftheta{\bfx\bfpi_m}}\!-\!1\!-\!\rmi\skal\bftheta {\bfx\bfpi_m} \eins_\DD(\bfx)\big)\,\XXX(\rmd \bfx)\,.\]
As $\skal{\bft\tr\bfmu}\bftheta=\sum_{m=1}^n \Delta t_{(m)}\skal{\bfmu\bfpi_m}\bftheta$
and $\|\bftheta\|^2_{\bft\tr\Sigma}=\sum_{m=1}^n \Delta t_{(m)}\|\bftheta\bfpi_m\|^2_\Sigma$, we get by recalling~\eqref{VVVsubord} and $\eins_\DD\circ \bfpi_m-\eins_\DD=(\eins_{\DD^C})(\eins_\DD\circ \bfpi_m)$ that
\begin{eqnarray*}
	\lefteqn{\sum_{m=1}^n\Delta t_{(m)}\int_{\RR_*^n}\left(\exp\{\rmi\skal\bftheta{\bfpi_m(\bfx)}\}\!-\!1\!-\!\rmi\skal\bftheta {\bfpi_m(\bfx)} \eins_\DD(\bfx)\right)\,\XXX(\rmd \bfx)}\nonumber&&\\
	&=&\rmi \skal{\bfc(\bft,\XXX)}\bftheta+\int_{\RR_*^n}\left(\exp\{\rmi\skal\bftheta{\bfx}\}\!-\!1\!-\!\rmi\skal\bftheta {\bfx} \eins_\DD(\bfx)\right)\,(\bft\tr\XXX)(\rmd \bfx)\,.
\end{eqnarray*}
By combining the above,~\eqref{multiexponent2} follows from~\eqref{multiexponent1}, completing the proof.\halmos
\subsection{Proof of Theorem~\ref{thmexist}}\label{subsecproofII}
We collect some useful estimates into a lemma. Its proof and purpose follow~\cite{s}~(see his Lemma 30.3) and~\cite{BPS01} (see the proof of their Theorem~3.2). However, we have to adapt these results to deal with the multivariate time parameter.
\begin{lemma}\label{lemcharlinbounf} If $X\sim L^n(\bfmu,\Sigma,\XXX)$ and $\bftheta\in\RR^n$, then there exist finite $C_1=C_1(\bftheta,X)$, $C_2=C_2(X)$ and
$C_3=C_3(X)$ such that, for $\bft\in[0,\infty)^n$,
\begin{eqnarray}
|\Phi_{X(\bft)}(\bftheta)-1|&\le& C_1(1\wedge\|\bft\|)\,,\label{eqsublin}\\
\EE[1\wedge \|X(\bft)\|^2]&\le& C_2 (1\wedge\|\bft\|)\,,\label{m2ineq}\\
\EE[1\wedge \|X(\bft)\|]&\le& C_2^{1/2} (1\wedge\|\bft\|^{1/2})\,,\label{m1ineq}\\
\big\|\EE[X(\bft)\eins_\DD(X(\bft))]\big\|&\le&C_3(1\wedge\|\bft\|)\label{m1truncabsoutineq}\,.
\end{eqnarray}
\end{lemma}
\noindent{\bf Proof.}~Let $\bft=(t_1,\dots,t_n)\in[0,\infty)^n$ and $\bftheta=(\theta_1,\dots,\theta_n)\in\RR^n$, and introduce a L\'evy measure $\NNN:=\sum_{\langle(1),\dots,(n)\rangle}\sum_{k=1}^n\XXX_{\{(k),\dots,(n)\}}$
with the first summation taken over all permutations $\langle(1),\dots,(n)\rangle$.

Recall 
$|e^z-1|\le|z|$, holds for $z\in\CC$ with $\Re z\le 0$, and, in particular, for $z:=\bft\tr\Psi(\bftheta)$ in~\eqref{multiexponent2}.
Further, we have $|\Re(\bft\tr\Psi(\bftheta))|\le C_{11}\|\bft\|$, where
\[C_{11}:=\frac 12\sum_{k,l=1}^n|\theta_k\theta_l\Sigma_{kl}|+\int_{\RR^n_*}\big|1\!-\!\cos\skal\bftheta\bfx\big|\;\NNN(\rmd\bfx)\,.\]
In~\eqref{defbfc}, note $\|\bfc(\bft,\XXX)\|\le n\XXX(\DD^C)\|\bft\|$, giving
$|\Im (\bft\tr\Psi(\bftheta))|\le C_{12}\|\bft\|$, where
\[C_{12}:=n(\|\bfmu\|+n\XXX(\DD^C))\|\bftheta\|+\int_{\RR^n_*}\big|\skal\bftheta\bfx \eins_\DD(\bfx)\!-\!\sin\skal\bftheta\bfx\big|\;\NNN(\rmd\bfx)\,.\]
Plainly,~$C_{11}$ and $C_{12}$ are finite constants in view of by~\eqref{Piintegrab}.
Choosing $C^2_{13}:=C^2_{13}(\bftheta):=C_{11}^2+C_{12}^2$ shows $|\Phi_{X(\bft)}(\bftheta)-1|\!\le\!C_{13}\|\bft\|$, so that~\eqref{eqsublin} holds for some finite $C_1=C_1(\bftheta)$.

Setting $\YYY_\bft(A):=(\bft\tr\XXX)(A\cap\!\DD^C)$ and $\ZZZ_\bft(A):=(\bft\tr\XXX)(A\cap\!\DD)$, $A\subseteq\RR^n_*$ Borel, yields L\'evy measures $\YYY_{\bft}$ and $\ZZZ_{\bft}$ on $\RR^n_*$ with disjoint supports and associated in\-de\-pen\-dent L\'evy pro\-cesses $Y^{(\bft)}\sim L^n({\bf 0},0,\YYY_\bft)$ and $Z^{(\bft)}=(Z^{(\bft)}_1,\dots,Z^{(\bft)}_n)\sim\!L^n(\bft\tr\bfmu+\bfc(\bft,\XXX),\bft\tr\Sigma,\ZZZ_\bft)$, respectively. By Proposition~\ref{propmultiparameter}, we may decompose $X(\bft)\eqd Y^{(\bft)}(1)+Z^{(\bft)}(1)$ into a sum of independent $n$-dimensional random vectors.

Note $Y^{(\bft)}$ is a compound Poisson process with jumps in $\|\cdot\|$-modulus larger than 1.
In particular,  $\{ Y^{(\bft)}\mbox{ has no jumps in time interval }[0,1]\}\subseteq\{Y^{(\bft)}(1)=0\}$,
giving the bound
\[\PP(Y^{(\bft)}(1)\neq 0)\le 1-\PP(Y^{(\bft)}\mbox{ has no jumps in time interval }[0,1]\}= 1-\exp(-(\bft\tr\XXX)(\DD^C))\,.\]
Since $(\bft\tr\XXX)(\DD^C)\le\|\bft\|\NNN(\DD^C)$ and $1-e^{-x}\le x$, $x\in\RR$, we have
\begin{equation}\label{Ybfnotnull}
\PP(Y^{(\bft)}(1)\neq 0)\,\le\,\NNN(\DD^C)\;\|\bft\|\,.
\end{equation}
On the other hand, $Z^{(\bft)}$ has jumps bounded in norm by 1. In particular, $Z^{(\bft)}(1)$
has finite moments of all order. Recall~$\EE[Z^{(\bft)}_k(1)]\!=\!\mu_kt_k\!+\!c_k(\bft,\XXX)$ and Var$(Z^{(\bft)}_k(1))\!=\!\Sigma_{kk}t_k\!+\!\int_{\DD_*}\! x^2_k(\bft\tr\XXX)(\rmd \bfx)$ for $1\!\le\! k\!\le\!n$ (see~\cite{s}, his Example~25.12).

By~\eqref{Piintegrab}, $C_{21}:=2\|\bfmu\|^2+2n^2 \XXX(\DD^C)^2+\spur(\Sigma)+\int_{\DD_*} \|\bfx\|^2\NNN(\rmd \bfx)\big\}$ is a finite constant, in addition satisfying
	\begin{equation}\label{msEuclid}
	\EE[\|Z^{(\bft)}(1)\|^2]
	\;\le\;C_{21}(\|\bft\|+\|\bft\|^2)\,.\end{equation}
By~\eqref{Ybfnotnull}--\eqref{msEuclid}, $\EE[1\wedge \|X(\bft)\|^2]\le C_{22} (\|\bft\|+\|\bft\|^2)$ holds with the choice $C_{22}:=\NNN(\DD^C)+C_{21}$, by noting $\EE[1\wedge \|X(\bft)\|^2]\le\PP(Y^{(\bft)}(1)\neq{\bf 0})+\EE[\|Z^{(\bft)}(1)\|^2]$. This completes the proof of~\eqref{m2ineq}, while~\eqref{m1ineq} is implied by~\eqref{m2ineq} and the Cauchy-Schwarz inequality.

Recall~$\|\bfz\|^2_\infty:=\max_{1\le k\le n}|z_k|^2\!\le\!\|\bfz\|\!:=\!\bfz\overline\bfz'$, $\bfz\!=\!(z_1,\dots,z_n)\!\in\!\CC^n$, and set $\DD_\infty\!:=\!\{\bfx\in\RR^n\!:\!\|\bfx\|_\infty\!\le\!1\}$. If $g(x):=e^{\rmi x}\!-\!1,x\in\RR$, we have
\begin{eqnarray*}
\lefteqn{\|\EE[X(\bft)\eins_{\DD_\infty}(X(\bft))]\|_{\infty} \le \max_{1\le j\le n}\big|\EE[g(X_j(\bft))\eins_{\DD^C_\infty}(X(\bft))]\big|}\\
&&+\displaystyle{\max_{1\le j\le n}}\big|\EE[(g(X_j(\bft))-\rmi X_j(\bft))\eins_{\DD_\infty}(X(\bft))]\big|+\displaystyle{\max_{1\le j\le n}}\big|\EE[g(X_j(\bft))]\big|\,.\end{eqnarray*}

By noting $\eins_{\DD^C_{\infty}}\!\le\!\eins_{\DD^C}\!\le\! 1\!\wedge\!\|\cdot\|^2$, we get
\[\big|\EE[g(X_j(\bft))\eins_{\DD^C_\infty}(X(\bft))]\big|\le2\EE[\eins_{\DD^C_\infty}(X(\bft))]\le2\EE[1\wedge \|X(\bft)\|^2]\,,\quad 1\le j\le n\,,\]
and then~\eqref{m2ineq} can be applied. Next, by noting~$4|g(x)\!-\!\rmi x|^2
\le x^4\!+\!x^6,x\in\RR$, we get
\[\big|\EE[(g(X_j(\bft))-\rmi X_j(\bft))\eins_{\DD_\infty}(X(\bft))]\big|
\,\le\,\EE[1\wedge X^2_j(\bft)]\,\le\,\EE[1\wedge \|X(\bft)\|^2]\,,\quad 1\le j\le n\,,
\]
and then~\eqref{m2ineq} can be applied. Lastly, we get $|\EE[g(X_j(\bft))]|\!=\!|\Phi_{X(\bft)}(\bfe_j)-1|, 1\!\le\!j\!\le\!n$, and then~\eqref{eqsublin} can be applied with $\bftheta\!\in\!\{\bfe_1,\dots,\bfe_n\}$.

Combining the above yields $\|\EE[X(\bft)\eins_{\DD_\infty}(X(\bft))]\|_{\infty}\le C_{31}(1\wedge\|\bft\|)$ for some finite constant $C_{31}$. Applying the Euclidean triangle inequality and $\|\cdot\|\le\!n^{1/2}\|\cdot\|_\infty$ yields
\[
\|\EE[X(\bft)\eins_{\DD}(X(\bft))]\|\le n^{1/2}\|\EE[X(\bft)\eins_{\DD_\infty}(X(\bft))]\|_{\infty}+\|\EE[X(\bft)\eins_{\DD_\infty\backslash\DD}(X(\bft))]\|\,.\]
The second term on the RHS is bounded from above by $n^{1/2}\EE[\eins_{\DD}(X(\bft))]$, and we found this to be bounded from above by $n^{1/2}\EE[1\!\wedge\! \|X(\bft)\|^2]$, to which~\eqref{m2ineq} was applicable. This completes the proof of the lemma.\halmos\\[1mm]
{\bf Proof of Theorem~\ref{thmexist}(i).}~Plainly, $\Theta$ in \eqref{SigmaY} is a valid covariance matrix as $\bfd\tr\Sigma$ is the covariance matrix of $B(\bfd)$ with $B\sim BM^n(0,\Sigma)$. It remains to validate that $\ZZZ$ in~\eqref{PliY}~is a L\'evy measure. By~\eqref{multiexponent2}, if $\bftheta\in\RR^n$, then~$\bft\mapsto (\bft\tr \Psi)(\bftheta)$ is a continuous function with domain $\bft\in[0,\infty)^n$. In particular, the family of probability measures $\{\PP(X(\bft)\in\rmd \bfx):\,\bft\in[0,\infty)^n\}$ is weakly continuous,
and $\PP(X(\bft)\in\rmd\bfx)$ is a Markov kernel from $[0,\infty)^n$ to $\RR^n$, and $\ZZZ_0(\rmd\bft,\rmd \bfx):=\eins_{[0,\infty)^n_*\times\RR^n}(\bft,\bfx)\PP(X(\bft)\!\in\!\rmd\bfx)\;\TTT(\rmd \bft)$ is a well-defined $\sigma$-finite Borel measure on the punctured product $([0,\infty)^n\times\RR^{n})_*$, for which we note
\begin{align}\int_{([0,\infty)^n\times\RR^{n})_*} \!\!1\!\wedge\!\|(\bft,\bfx)\|^2\,\ZZZ_0(\rmd\bft,\rmd \bfx)\!=\!
	\int_{[0,\infty)^n_*}\!\EE[1\!\wedge\!\|(I\bfe,X)(\bft)\|^2]\,\TTT(\rmd\bft)\,,\label{ineqtildePiLevmeasure}
\end{align}
where $\bfe=(1,\dots,1)\in\RR^n$ and $(I\bfe,X)$ is an augmented $2n$-dimen\-sional L\'evy process.

For $\bft\in[0,\infty)^n$, by noting $\|(\bft,\bft)\|^2=2\|\bft\|^2$ and $1\wedge \|(\bft,\bft)\|\!\le\!2^{1/2}\!(1\wedge\|\bft\|),\bft\in[0,\infty)^n$, and applying~\eqref{m2ineq}
with $C_2:=C_2((I\bfe,X))$, we get
\begin{eqnarray*}
	\EE[1\wedge \|(I\bfe,X)(\bft)\|^2]\,\le\,
C_2\,(1\wedge \|(\bft,\bft)\|)\le 2^{1/2}C_2\,(1\wedge \|\bft\|)\,.\end{eqnarray*}
As~\eqref{PiFVintegrab} holds for $\TTT$, the RHS in the last display is~$\TTT$-integrable, hence $\ZZZ_0$~and $\ZZZ$ in~\eqref{PliY} are L\'evy measures by~\eqref{ineqtildePiLevmeasure}.

Note $\|\bft\|\,\PP((\bft,X(\bft))\in\DD)\le \|\bft\|\,\eins_\DD(\bft)$
for all $\bft\in[0,\infty)^n$. As the RHS is $\TTT$-integrable by~\eqref{PiFVintegrab}, so is the LHS, and then~\eqref{defmu1} is well-defined.
The RHS of~\eqref{defmu2} is well-defined as an implication of~\eqref{m1truncabsoutineq}, applied to the augmented process~$(I\bfe,X)$.\\
{\bf Proof of Theorem~\ref{thmexist}(ii).}~On a suitable augmentation of $(\Omega,\FFF,\PP)$, where $T$ lives, we find $W\sim L^{2n}(\bfm,\Theta,\bfdelta_{{\bf 0}}\otimes(\bfd\tr \XXX))$, $\bfm=(\bfm_1,\bfm_2)$ with $\bfm_1$, $\bfm_2$ and $\Theta$ as in \eqref{defmu1}-- \eqref{SigmaY}, and a family $\xi=\{\xi(t,\bft):(t,\bft)\in[0,\infty)\times[0,\infty)^n_*\}$ of independent random vectors, satisfying $\xi(t,\bft)\eqd X(\bft)$
for $(t,\bft)\in[0,\infty)\times[0,\infty)^n_*$, such that $T,\xi,W$ are independent. Introduce a marked Poisson point process \[\ZZ_0:=\sum_{t>0}\bfdelta_{(t,T(t)-T(t-),\xi(t,T(t)-T(t-)))}\,,\]
thus being a Poisson point process with intensity $\rmd t \otimes \ZZZ_0$, where $\ZZZ_0$ is the L\'evy measure in Part~(i). Particularly,~$\ZZ_0$ is the point measure of jumps of a L\'evy process $Z_0\sim L^{2n}({\bf 0},{\bf 0},\ZZZ_0)$
via its L\'evy-It\^o decomposition. 
As $Z_1= T$, $Z=(Z_1,Z_2):=Z_0+W\eqd (T,X\odot T)$ is $T$ subordinating $X$ in the semi-strong sense.\\
{\bf Proof of Theorem~\ref{thmexist}(iii).}~If, in addition, $\int_{[0,1]_*^n}\|\bft\|^{1/2}\;\TTT(\rmd\bft)$ is finite, then~\eqref{PiFVintegrab} holds.
This follows similarly as in~\eqref{ineqtildePiLevmeasure}, but using~\eqref{m1ineq} instead of~\eqref{m2ineq}.
\halmos
\subsection{Proof of Propositions~\ref{propindepindist} and~\ref{propmonotindepmarginal}}\label{subsecproofIII}
Let $T\!=\!(T_1,T_2)\!\sim\! S^2$ and $X\!=\!(X_1,X_2)\!\sim\! L^2$ be independent. For $\bftheta=(\theta_1,\theta_2)\in\RR^2$, $t\!\ge\! s\!\ge\!0$, introduce $\widehat\Psi_X(\bftheta):=
\Psi_X(\bftheta)-\Psi_{X_1}(\theta_1)-\Psi_{X_2}(\theta_2)$.\\[1mm]
\noindent{\bf Proof of Propositions~\ref{propindepindist}.}~For $\bftheta=(\theta_1,\theta_2)\in\RR^2$, $t\!\ge\! s\!\ge\!0$, introduce $A(s,t)\!:=\!(T_1(s)\!\wedge\! T_2(t))\!-\!(T_1(s)\!\wedge\! T_2(s))$ and
$Z(s,t,\bftheta):=T_1(s)\Psi_{X_1}(\theta_1)+(T_2(t)\!-\!T_2(s))\Psi_{X_2}(\theta_2)$.

In view of~\eqref{multiexponent1}, for $\bftheta=(\theta_1,\theta_2)\in\RR^2$, $r\!\ge\! 0$, $t\!\ge\! s\!\ge\!0$, note
\[
(r,t,s)\tr \Psi_{X_1,X_2,X_2}(\bftheta,-\theta_2)=r\Psi_{X_1}(\theta_1)+(t\!-\!s)\Psi_{X_2}(\theta_2)+\widehat\Psi_X(\bftheta)(r\wedge t-r\wedge s)\,,
\]
and thus, by conditioning on $T$,
\[\Phi_{(X_1(T_1(s)),X_2(T_2(t))\!-\!X_2(T_2(s)))}(\bftheta)
=\EE\exp\{Z(s,t,\bftheta)+\widehat\Psi_X(\bftheta) A(s,t)\}\,.\]
As $X\circ T$ is assumed to be a L\'evy process, both $T$ and $X\circ T$, have independent
increments across the components. Conditioning the LHS of the last display on $T$ shows the following identity, for $\bftheta=(\theta_1,\theta_2)\in\RR^2,t\!\ge\! s\!\ge\!0$,
\begin{equation}\label{notleevy5}
\EE\exp\{Z(s,t,\bftheta)\}\,=\,\EE\exp\{Z(s,t,\bftheta)+\widehat\Psi_X(\bftheta) A(s,t)\}\,.
\end{equation}
(i)~Assume $X\eqd-X$. Since $X_1$ and $X_2$ are dependent, there exist $\bftheta=(\theta_1,\theta_2)\in\RR^2$ such that
$\widehat\Psi_X(\bftheta)\neq 0$. By symmetry, $\Psi_X(\bftheta),\widehat\Psi_X(\bftheta),\Psi_{X_j}(\theta_j)\in\RR$, $j=1,2$.
Let $t\!>\!0, u\!\ge\!1$. In \eqref{notleevy5} we have $Z(t,ut,\bftheta)\in\RR$, forcing $A(t,ut)=0$ almost surely.
In particular, $u\mapsto A(t,ut)$ degenerates to the null process. As $T_2$ cannot degenerate to the null process, we must have
$T_2(t)<T_2(ut)$ for some $u>1$ with probability one, and thus, $T_1(t)\le T_2(t)$ almost surely. Reversing the role of $T_1$ and $T_2$ completes the proof of Part~(i).\\
(ii)~As $X_1,X_2$ are dependent we have $\widehat\Psi_X(\bftheta)\neq 0$ for some $\bftheta\in\RR^2$. If $T$ is deterministic with drift $(d_1,d_2)$, then $\widehat\Psi_X(\bftheta)A(t,(1\!+\!\varepsilon)t)\in2\pi\rmi\ZZ$ for $t,\varepsilon>0$,
as an implication of~\eqref{notleevy5}, giving $d_{1}\le d_{2}$, with the argument being completed as in (i).\\
(iii)~Assume $T(1)$, and thus $A(t,2t)$ for all $t\geq0$, admits a finite first moment. In addition, suppose there exists a sequence $\bftheta_n\to {\bf 0}$ as $n\!\to\!\infty$ such that
$\widehat \Psi(\bftheta_n)\neq 0$ and $\Re\widehat \Psi(\bftheta_n)\le 0$ for all $n$. As $|1-e^z|\le |z|$ for $\Re z\le 0$, note $|Z(t,2t,\bftheta_n)(1-\exp\{A(t,2t)\widehat\Psi_X(\bftheta_n)\})/\widehat \Psi(\bftheta_n)|\le A(t,2t)$, and dominated convergence is applicable to \eqref{notleevy5}, giving $A(t,2t)=0$ almost surely, since
\[
0=\lim_{n\to\infty}\EE[Z(t,2t,\bftheta_n)(1\!-\!\exp\{A(t,2t)\widehat\Psi_X(\bftheta_n)\})/\widehat \Psi_X(\bftheta_n)\big]=\EE[A(t,2t)]\,.
\]
If $X(1)$ admits a finite second moment, then $\widehat \Psi_X(\bftheta)=-\rho \theta_1\theta_2 +o(\|\bftheta\|^2)$ as $\bftheta\to {\bf 0}$, where $\rho=$Cov$(X_1(1),X_2(1))$, the existence of sequence as required in the previous paragraph is obvious, provided $\rho\neq 0$.\halmos\\[1mm]
\noindent{\bf Proof of Proposition~\ref{propmonotindepmarginal}.}~Introduce $D:=T_2-T_1\sim FV^1(d,\DDD)$.
If $\bftheta=(\theta_1,\theta_2)\in\RR^2$, note $(r,s)\tr \Psi_{X}(\bftheta)=(r\!\wedge\!s)\Psi_{X}(\bftheta)\!+\!(s\!-\!r)^+\Psi_{X_2}(\theta_2)\!+\!(s\!-\!r)^-\Psi_{X_1}(\theta_1),r,s\!\ge\!0$, so that, by conditioning on $T$,
\begin{equation}
\Phi_{X(T(t))}(\bftheta)
=\EE[\exp\{(T_1(t)\!\wedge\!T_2(t))\Psi_{X}(\bftheta)\!+\!D^+(t)\Psi_{X_2}(\theta_2)\!+\!D^-(t)\Psi_{X_1}(\theta_1)\}]\,.\label{eqPhiXT}\end{equation}
(i)~Recall $T$ is monotonic if and only if either $D$ or $-D$ is a subordinator. As we assumed $T$ to have non-monotonic and non-deterministic components,
one of the following exclusive cases holds~(see \cite{s}, his Corollary~24.8 and his Theorem~24.10):\\
(a)~$\DDD(-\infty,0)>0$, $\DDD(0,\infty)=0$ and $d>0$, so that the support of the distribution of $D(1)$ is unbounded towards $-\infty$ with $d$ as its supremum;\\
(b)~$\DDD(-\infty,0)=0$, $\DDD(0,\infty)>0$ and $d<0$, so that the support of $D(1)$ is unbounded towards $\infty$ with $d$ as its infimum;\\
(c)~$\DDD(-\infty,0)>0$ and $\DDD(0,\infty)>0$ ($d\in\RR$ is arbitrary), so that the support of $D(1)$ is unbounded towards $\infty$ and $-\infty$.\\[1mm]
In all cases, we have $\PP(D(1)\!>\!0)\!>\!0$ and $\PP(D(1)\!<\!0)\!>\!0$, implying $\EE[D^+(1)]\!>\!0$ and $\EE[D^-(1)]\!>\!0$, respectively. We assumed a finite
second moment for $T$, so that the second moment of $D$ is finite, implying $\EE[|D(1)|]=\EE[D^+(1)]\!+\!\EE[D^-(1)]<\infty$.

Assume $\EE[D^+(t)]=t\EE[D^+(1)]$, $t\ge 0$, so that
$\EE[D^{+}(1)]=\EE[(D(n)/n)^{+}]$, $n=1,2,3,\dots$. Consequently, we have
$\lim_{n\to\infty}\EE[(D(n)/n)^+]=\EE[\EE[D(1)]^{+}]=\EE[D(1)]^{+}$, as convergence in mean holds in the context of the
strong law of large numbers for independent and identically distributed integrable random variables. This leads to the contradiction
$\EE[D^{+}(1)]=\EE[D(1)]^+=(\EE[D^+(1)]-\EE[D^-(1)])^+<\EE[D^+(1)]$. To summarise, $t\mapsto\EE[D^+(t)],t\ge 0$ cannot be a linear function.

On the RHS of~\eqref{eqPhiXT}, taking partial derivatives twice
with respect to $\bftheta=(\theta_1,\theta_2)$ 
under the expectation and applying dominated convergence to $\bftheta\to{\bf 0}$, 
we derive the Wald-type identity
\begin{equation}\label{covstppedXT}
\myCov(X_1(T_1(t)),X_2(T_2(t)))=
\EE[X_1(1)]\EE[X_2(1)]\,\myCov(T_1(t),T_2(t))+\rho\EE[T_1(t)\wedge T_2(t)]\,.
\end{equation}
By our assumptions, $T$ and $X$ admit finite second moments, so that both sides of~\eqref{covstppedXT} are finite.

Contradicting the hypothesis, assume $X(T(t))\!\eqd\!Y(t)$, for all $t\!\ge0$, where $Y$ is a given bivariate L\'evy process. Plainly, $T$ and $Y$ are L\'evy processes with finite second moments. In particular, $t\mapsto\myCov(T_1(t),T_2(t))$ and $t\mapsto\myCov(Y_1(t),Y_2(t))$
are linear functions, and so is $t\mapsto\EE[T_1(t)\wedge T_2(t)]$, as we assumed $\rho\!\neq\!0$ in~\eqref{covstppedXT}.

Also, $t\mapsto\EE[T_2(t)]$ is linear, so that noting $\EE[T_1(t)\wedge T_2(t)]=
\EE[T_2(t)]-\EE[D^+(t)]$, $t\ge0$, contradicts the non-linearity of $t\mapsto\EE[D^+(t)]$, completing the proof of~(i).\\
(ii)~If $T_1,T_2$ are independent and drift-less, the components of
$X\odot\!T$ are independent by Proposition~\ref{propsubordindependent}. Then using Proposition~\ref{propsupextendssub} on each component yields
$X\odot T\eqd (X_1\circ\!T_1,X_2^*\circ\!T_2)$ for independent L\'evy processes $T_1,T_2,X_1,X_2^*$, where $X_2^*\eqd X_2$.
To summarise, we have 
\begin{equation}\label{Phiprod}\Phi_{X\odot T(t)}(\bftheta)=\EE[\exp\{T_1(t)\Psi_{X_1}(\theta_1)+T_2(t)\Psi_{X_2}(\theta_2)\}],\quad
\bftheta=(\theta_1,\theta_2)\in\RR^2\,.\end{equation}
Next, note $(r\wedge s)z+(s-r)^+z_2+(s-r)^-z_1=(r\wedge s)\widehat z+(rz_1+sz_2), r,s\ge 0,z,z_1,z_2,\widehat z:=z\!-\!z_1\!-\!z_2\in\CC$.

As we assume that $X_1,X_2$ are dependent, there exists $\bftheta^*\in\RR^2$ such that $\widehat\Psi_X(\bftheta^*)\neq 0$.
Further, $\Psi_{X_j}(\theta^*_j)\le 0, j=1,2,\widehat\Psi_{X}(\bftheta^*)\in\RR$ by our symmetry assumption $X\eqd-X$. If for all $t>0$,~\eqref{Phiprod} matches~\eqref{eqPhiXT}, we have \[\EE[\exp\{T_1(t)\Psi_{X_1}(\theta^*_1)+T_2(t)\Psi_{X_2}(\theta^*_2)\}(\exp\{(T_1(t)\wedge T_2(t))\widehat\Psi_{X}(\bftheta^*)\}-1)]=0\,,\]
with the implication $T_1(t)\wedge T_2(t)=0$, a.s., for all $t>0$. In particular,
the null process and $T_1\wedge T_2$ must be indistinguishable as processes, which is a contradiction to $T_1,T_2$ being nontrivial subordinators, completing the proof.\halmos
\subsection{Proof of Proposition~\ref{propVGGCFV}}\label{subsecproofIV}
(i)~See \cite{BKMS16}, Part (a) of their Proposition 2.1.\\[1mm]
(ii)~Let~$0\!<\!q\!<\!1$ and~$B\!\sim\! BM^n(\bfmu,\Sigma)$. If $\bft\!\in\![0,\infty)^n_*$, set~$\psi(\bft)\!:=\!\EE[\|B(\bft)\|^{2q}\eins_{\DD_*}(B(\bft))]$.\\
`$\Rightarrow$': Note $\psi(\bft)\le 1$ for $\bft\!\in\![0,\infty)^n$, and introduce $|\Sigma|:=(|\Sigma_{kl}|)\in\RR^{n\times n}$, $|Z|:=(|Z_1|,\dots,|Z_n|)$ for a standard normal vector $(Z_1,\dots,Z_n)\in\RR^n$
and $Q:=2^{(2q)\vee 1}\,(\|\bfmu\|^{2q}+ \EE[\|\,|Z|\,\|^{2q}_{|\Sigma|}]\big)$.
If $\bft\in [0,\infty)^n_*$, we have
$\psi(\bft)\le \EE[\|B(\bft)\|^{2q}]\le Q(\|\bft\|^{q}\vee\|\bft\|^{2q})$
and thus $\psi(\bft)\!\le\! (1\vee Q) (1\wedge \|\bft\|^q)$. As~$\YYY(\rmd \bfy)\allowbreak= \PP(B(\bft)\in \rmd \bfy)\,\TTT(\rmd \bft)$ in \eqref{PliY2},~sufficiency follows from this and Part~(i).\\
`$\Leftarrow$':~Assume an invertible $\Sigma$. Set $\DD^+_*:=\DD\cap[0,\infty)^n_*$.
The proof is completed, provided
we can show that $i:=\inf_{\bft\in\DD^+_*} \psi(\bft)/\|\bft\|^{q}>0$. If $\bft\in[0,\infty)_*^n$, then set $\phi(\bft):=\psi(\bft)/\|\bft\|^{q}$.
Plainly, we have
$\bft_m\to\bft_0$, $\bfs_m:=\bft_m/\|\bft_m\|\to\bfs_0$ and $\phi(\bft_m)\to i$ as $m\to\infty$
for some $\bft_0\in\DD,\bfs_0\in \mySS_+,\bft_m\in\DD_*,m\ge 1$.

If $\bft_0\neq{\bf 0}$, then we find $\emptyset\neq J\subseteq \{1,\dots,n\}$ such that $\bft_0\in C_J:=\{\sum_{j\in J} x_j \bfe_j:x_j\!>\!0\;\mbox{for all }j\!\in\!J\}$. Note
$\bft_0\tr\bfmu \in \bfpi_J(\RR^n)$, while
$\bft_0\tr\Sigma:\bfpi_J(\RR^n)\to\bfpi_J(\RR^n)$ is invertible. Particularly, $\PP(B(\bft)\neq {\bf 0})=1$,
$\PP(0\!<\!\|B(\bft)\|\!<\!1)\!>\!0$ and $\psi^*(\bft):=\EE[\|B(\bft)\|^{2q}\eins_{(0,1)}(\|B(\bft)\|)]\allowbreak>0$.
As desired, we get from Fatou's lemma and the continuity of the sample paths of $B$ that $i=\|\bft_0\|^{-q}\liminf_{m\to\infty}\psi(\bft_m)\ge \|\bft_0\|^{-q}\psi^*(\bft_0)>0$.

If $\bft_0={\bf 0}$, let $B^*\!\sim\!BM^n({\bf 0},\Sigma)$, and recall $\|\bft_m\|^{-\frac{1}{2}} B(\bft_m)\eqd \bfmu\tr(\|\bft_m\|^{\frac{1}{2}}\bfs_m)+B^*(\bfs_m)=:W_m\to B^*(\bfs_0)$ and $\eins_{\DD_*}(\|\bft\|^{1/2} W_m)=\eins_{\DD}(\|\bft\|^{1/2} W_m)$$\to 1$, almost surely, as $m\to\infty$, by continuity of the sample paths of $B^*$. The proof of the necessity is completed by Fatou's lemma as \[i\ge \liminf_{m\to\infty} \EE[\|W_m\|^{2q}\eins_{\DD_*}(\|\bft\|^{1/2} W_m)]\ge \EE[\|B^*(\bfs_0)\|^{2q}]>0\,.\]
\noindent(iii)~Let $Y\sim \VGGC^n(\bfd,\bfmu,\Sigma,\UUU)$. To have $Y\sim FV^n$ for an invertible $\Sigma$, we cannot allow for a non-trivial Brownian component in~\eqref{GVGcharexpo}, thus forcing $\bfd={\bf 0}$. The L\'evy measure of an  $FV^n$-process obeys~\eqref{PiFVintegrab}, and the necessity part of Proposition~\ref{propVGGCFV}(ii)
forces~$\int_{\DD^C}\UUU(\rmd\bfu)/\|\bfu\|^{1/2}$ to be finite for invertible $\Sigma$.
If $\bfd={\bf 0}$ then $Y$  in~\eqref{GVGcharexpo} has no Brownian component. If $\int_{\DD^C}\UUU(\rmd\bfu)/\|\bfu\|^{1/2}$ is finite, then so is~\eqref{PiFVintegrab}
for $\YYY$, as an implication of Proposition~\ref{propVGGCFV} for $q=1$.
If $\bfd={\bf 0}$ and~$\int_{\DD^C}\UUU(\rmd\bfu)/\|\bfu\|^{1/2}<\infty$ hold simultaneously with no invertibility assumptions on $\Sigma$, then so does $Y\sim FV^n$.\halmos

\end{document}